%% file: block_smoothers.tex
\documentclass[a4paper, 10pt]{article}

\usepackage[utf8]{inputenc}

\usepackage{amsmath,amssymb,amsfonts}
\usepackage{mathtools}
\usepackage{booktabs}
\usepackage{bm} % bold math

\usepackage[usenames,dvipsnames]{xcolor}
\usepackage{graphicx}
\usepackage{blindtext}

\usepackage{siunitx}
\DeclareSIUnit\flops{FLOP\per\second}
\DeclareSIUnit\ups{Up\per\second}
\DeclareSIUnit\lups{LUP\per\second}
\DeclareSIUnit\lu{LU}
\DeclareSIUnit\byte{B}
\DeclareSIUnit\bit{b}
\DeclareSIUnit\doubles{DP}
\DeclareSIUnit\cycles{cy}

\usepackage{todonotes}
\usepackage[labelformat=simple]{subcaption}

\usepackage{tikz}
\usetikzlibrary{calc,3d}
\usepackage{calc,pgf,pgfplots}
\usepackage{tikzscale}

\usepackage{listings}

\usepackage{calc}
\usepackage{algorithm}
\usepackage{algpseudocode}
\let\oldReturn\Return
\renewcommand{\Return}{\State\oldReturn}

\usepackage{setspace}

\renewcommand{\Comment}[2][.5\linewidth]{%
  \leavevmode\hfill\makebox[#1][l]{$\triangleright$~#2}}

\usepackage{prettyref}
\newrefformat{tab}{Table~\ref{#1}}
\newrefformat{fig}{Figure~\ref{#1}}
\newrefformat{alg}{Algorithm~\ref{#1}}
\newrefformat{sec}{Section~\ref{#1}}

\usepackage{geometry}
\lstdefinestyle{common-style}{
	basicstyle=\footnotesize,       % the size of the fonts that are used for the code
	numbers=left,                   % where to put the line-numbers
	numberstyle=\footnotesize,      % the size of the fonts that are used for the line-numbers
	stepnumber=1,                   % the step between two line-numbers. If it is 1 each line will be numbered
	numbersep=5pt,                  % how far the line-numbers are from the code
	showspaces=false,               % show spaces adding particular underscores
	showstringspaces=false,         % underline spaces within strings
	showtabs=false,                 % show tabs within strings adding particular underscores
	frame=single,                   % adds a frame around the code
	tabsize=2,                      % sets default tabsize to 2 spaces
	captionpos=b,                   % sets the caption-position to bottom
	breaklines=false,                % sets automatic line breaking
	breakatwhitespace=false,        % sets if automatic breaks should only happen at whitespace
	keywordstyle={\color{blue}\textbf},         % keywords are blue, (and blue)
	commentstyle={\color{OliveGreen}},          % comments
	literate={\$}{{\$}}1,
	escapechar=\&
}

\lstdefinestyle{fortran-style}{
	morecomment=[l][\color{\colorpragmaacc}]{!\$acc}, % ACC directives are orange
	morecomment=[l][\color{\colorpragmaomp}]{!\$omp},          % OMP directives are red
}

\usepackage{authblk}
\graphicspath{{.} {build/}}

\makeatletter
\def\title#1{\gdef\@title{#1}\gdef\THETITLE{#1}}
\makeatother
\input{commands}

\author[a,b]{Immo Huismann\thanks{Corresponding author: Immo.Huismann@tu-dresden.de}}
\author[a,b]{Jörg Stiller}
\author[a,b]{Jochen Fröhlich}
\affil[a]{Institute of Fluid Mechanics, TU Dresden}
\affil[b]{Center for Advancing Electronics Dresden (cfaed)}
\title{Scaling to the stars -- a linearly scaling elliptic solver for~$\poly$-multigrid}

\begin{document}

\maketitle

\begin{abstract}
  High-order methods gain increased attention in computational fluid dynamics.
  However, due to the time step restrictions arising from the semi-implicit time stepping for the incompressible case, the potential advantage of these methods depends critically on efficient elliptic solvers.
  Due to the operation counts of operators scaling with with the polynomial degree~$\poly$ times the number of degrees of freedom~$\ndof$, the runtime of the best available multigrid solvers scales with~$\order{\poly \cdot \ndof}$.
  This scaling with~$\poly$ significantly lowers the applicability of high-order methods to high orders.
  While the operators for residual evaluation can be linearized when using static condensation, \schwarz-type smoothers require their inverses on fixed subdomains.
  No explicit inverse is known in the condensed case and matrix-matrix multiplications scale with~${\poly \cdot \ndof}$.
  This paper derives a matrix-free explicit inverse for the static condensed operator in a cuboidal subdomain.
  It scales with~$\poly^3$ per element, i.e.\,${\ndof}$ globally, and allows for a linearly scaling additive~\schwarz\ smoother, yielding a $\poly$-multigrid cycle with an operation count of~$\order{\ndof}$.
  The resulting solver uses fewer than four iterations for all polynomial degrees to reduce the residual by ten orders and has a runtime scaling linearly with~$\ndof$ for polynomial degrees at least up to~$48$.
  Furthermore the runtime is less than one microsecond per unknown over wide parameter ranges when using one core of a CPU, leading to time-stepping for the incompressible \navierstokes\ equations using as much time for explicitly treated convection terms as for the elliptic solvers.
\end{abstract}

\input{block_smoothers_introduction}
% \newpage
\input{block_smoothers_discretization}
% \newpage
\input{block_smoothers_smoothers}
% % \newpage
% \input{block_smoothers_linear_inverse}
% \newpage
\input{block_smoothers_multigrid}
% \newpage
\input{block_smoothers_operator_runtimes}
\input{block_smoothers_solver_tests}

% \newpage
\input{block_smoothers_extensions}

% \newpage
\input{block_smoothers_conclusion}

\subsubsection*{Acknowledgements:}
This work is supported in part by the German Research Foundation (DFG) within the Cluster of Excellence ‘Center for Advancing Electronics Dresden’ (cfaed).
The authors would like to thank their colleagues in the Orchestration path of cfaed for stimulating discussions and ZIH, Dresden, for the provided computational resources and M.\,Kronbichler who kindly provided the reference data for the~\taylorgreen\ vortex.

\bibliographystyle{abbrv}
\bibliography{block_smoothers.bib}

\end{document}

%% file: commands.tex
\newcommand{\of}[1]{\!\left(#1\right)}

\newcommand{\twall}{t_{\mathrm{wall}}}
\newcommand{\name}[1]{\textsc{#1}}
\newcommand{\helmholtz}{\name{Helmholtz}}
\newcommand{\laplace}{\name{Laplace}}
\newcommand{\poisson}{\name{Poisson}}
\newcommand{\dirichlet}{\name{Dirichlet}}
\newcommand{\neumann}{\name{Neumann}}

\newcommand{\taylorgreen}{\name{Taylor}-\name{Green}}
\newcommand{\moura}{\name{Moura}}
\newcommand{\schwarz}{\name{Schwarz}}
\newcommand{\schur}{\name{Schur}}
\newcommand{\galerkin}{\name{Galerkin}}
\newcommand{\krylov}{\name{Krylov}}
\newcommand{\navierstokes}{\name{Navier-Stokes}}
\newcommand{\reynolds}{\name{Reynolds}}

\newcommand{\GLL}{GLL}
\newcommand{\specht}{SPECHT\_FS}
\newcommand{\scorep}{\mbox{Score-P}}

\newcommand{\domain}{\Omega}

\newcommand{\helmholtzopsymb}{H}
\newcommand{\helmholtzop}{\mathbf{\helmholtzopsymb}}
\newcommand{\laplaceop}{\Delta}

\newcommand{\weighsymbol}{W}
\newcommand{\restrictsymbol}{R}

\newcommand{\solverdtcg}{dCG}
\newcommand{\solvermg}{MG}
\newcommand{\solvermgcg}{kMG}
\newcommand{\solverkcvmg}{kvMG}

\newcommand{\ntime}{n_{\mathrm{t}}}

\newcommand{\elementindex}{e}
\newcommand{\nelement}{n_{\mathrm{\elementindex}}}
\newcommand{\nvertex}{n_{\mathrm{v}}}
\newcommand{\ndof}{n_{\mathrm{DOF}}}
\newcommand{\ndofstar}{n_{\mathrm{DOF}}^{\star}}
\newcommand{\ncores}{n_{\mathrm{cores}}}

\newcommand{\poly}{p}

\newcommand{\eqdot}{\quad .}
\newcommand{\eqcomma}{\quad ,}
\newcommand{\varu}{u}

\newcommand{\varf}{f}

\newcommand{\mymat}[1]{\mathbf{#1}}

\newcommand{\vecu}{\mymat{u}}

\newcommand{\vecF}{\mymat{F}}

\newcommand{\bound}{\mathrm{B}}
\newcommand{\inner}{\mathrm{I}}
\newcommand{\face}{\mathrm{F}}

\newcommand{\transmatsymb}{S}
\newcommand{\transmat}{\mymat{\transmatsymb}}

\newcommand{\masssymb}{M}
\newcommand{\stiffsymb}{L}
\newcommand{\standardsymb}{\mathrm{S}}

\newcommand{\massmat}{\mymat{\masssymb}}
\newcommand{\stiffmat}{\mymat{\stiffsymb}}

\newcommand{\standardmassmat}{\mymat{\masssymb}_{\standardsymb}}
\newcommand{\standardstiffmat}{\mymat{\stiffsymb}_{\standardsymb}}

\newcommand{\nstar}{n_{\mathrm{S}}}

\newcommand{\identitymat}{\mymat{I}}
\newcommand{\eigenvaluemat}{\mymat{\Lambda}}

\newcommand{\tptwo}[2]{#1 \otimes #2}

\newcommand{\tp}[3]{#1 \otimes #2 \otimes #3}

\newcommand{\ptp}[3]{\left(#1 \otimes #2 \otimes #3\right)}
\newcommand{\geomcoeffsymb}{d}
\newcommand{\geomcoeff}[1]{\geomcoeffsymb_{#1}}

\newcommand{\diagonal}{\mymat{D}}

\newcommand{\order}[1]{\mathcal{O}\!\left(#1\right)}
\newcommand{\eigenspace}{\mathrm{E}}

\newcommand{\condop}{\mymat{\hat{\helmholtzopsymb}}}

\newcommand{\condrestrict}{\mymat{\hat{\restrictsymbol}}}
\newcommand{\condweigh}{\mymat{\hat{\weighsymbol}}}
\newcommand{\condinterpol}{\bm{\hat{\mathcal{J}}}}

\newcommand{\condrhs}{\mymat{\hat{F}}}
\newcommand{\condvar}{\mymat{\hat{u}}}
\newcommand{\starrhs}{\mymat{\tilde{F}}}
\newcommand{\starvar}{\mymat{\tilde{u}}}

\newcommand{\condF}{\mymat{\hat{F}}}
\newcommand{\condu}{\mymat{\hat{u}}}
\newcommand{\condp}{\mymat{\hat{p}}}
\newcommand{\condq}{\mymat{\hat{q}}}
\newcommand{\condr}{\mymat{\hat{r}}}
\newcommand{\conds}{\mymat{\hat{s}}}
\newcommand{\condz}{\mymat{\hat{z}}}

\newcommand{\weigh}{\mymat{\weighsymbol}}

\newcommand{\istar}{\mathrm{S}}

%% file: block_smoothers_introduction.tex
\section{Introduction}%
\label{sec:introduction}

Current focus in Computational Fluid Dynamics~(CFD) is on high-order methods combining the geometric flexibility of low-order methods such as the Finite Volume method with the stellar convergence properties of traditional Fourier methods.
%
% While the Discontinuous Galerkin method is receiving the most attention at the moment~\cite{hesthaven_2007_dg}, the spectral-element method~(SEM)~\cite{deville_2002_sem, karniadakis_1999_sem} and $h/p$-FEM have the same convergence properties.
% \todo{so what?}
%
Continuous as well as discontinuous~\galerkin\ methods using spectral or $h/p$-elements are the most prominent members of this group.
While for low-order methods the error scales with the element width~$h$ squared, it scales with~$h^{\poly}$ with high-order methods.
This spectral convergence property fuels the race to ever higher polynomial degrees.
Where the first application of the spectral-element method~(SEM) utilized~${\poly=6}$~\cite{patera_1984_sem}, typical polynomial orders in current simulations range up to~$11$~\cite{atak_2016_sim, beck_2014_dg,  lombard_2015_sim, merzari_2013_sim, serson_2017_cfd} and even~$p=15$ is not uncommon~\cite{merzari_2013_sim,fehn_2018_sim}.
However, while the higher polynomial degrees lead to faster convergence, the runtime of operators and, hence, solvers scale super-linearly with the polynomial degree, even with the exploitation of tensor-product bases.

For the simulation of incompressible fluid flow the solution of the elliptic~\helmholtz\ equation~${\lambda u - \Delta u = f}$ is required multiple times per time step: three times for the diffusion terms with positive~$\lambda$, and once for the pressure involving the harder case~${\lambda = 0}$.
For low-order methods, multigrid has been established as a very efficient solution technique.
In high-order methods, $\poly$-multigrid, where the polynomial degree is lowered and raised, offers an attractive alternative to~$h$-multigrid.
In both cases, smoothers are required and overlapping~\schwarz\ smoothers are a good option~\cite{haupt_2013_mg, lottes_2005_mg, stiller_2017_multigrid}.
These decompose the grid into overlapping blocks on which the exact inverse is applied, producing exceptional smoothing rates.
But these methods not only require the super-linearly scaling operator for residual evaluation, but a super-linearly scaling solver as well.
This super-linear scaling prevents the respective multigrid techniques from achieving their full potential and prohibits the usage of large polynomial degrees.

The goal of this paper is to devise a linearly scaling multigrid cycle for three-dimensional structured Cartesian grids, where residual evaluation, smoothing, restriction, and prolongation all scale linearly with the number degrees of freedom, without a further factor of~$\poly$ present in other methods.
% Second, they need to be linearized such that they scale with the
%
To this end the residual evaluation derived in~\cite{huismann_2017_condensation} for the static condensed case is combined with the~$p$-multigrid block smoother technique proposed in~\cite{haupt_2013_mg, haupt_2017_mg}, leaving only the inverse on a~$2^3$ element block as non-matrix-free, super-linearly scaling operator.
This operator is embedded in the full system and then factorized and rearranged, attaining a matrix-free inverse that is then factorized to linear complexity.
The inverse leads to a linearly scaling multigrid cycle, which is confirmed in runtime tests ranging from solution of the~\poisson\ equation to benchmarks for the~\navierstokes\ equations.

The paper is structured as follows: First, the main steps of the discretization are recalled, then in~\prettyref{sec:smoothers}, the linearly scaling inverse of the local~\schwarz\ operators is derived.
Afterwards runtime tests for the smoothers confirm the linear scaling and the resulting multigrid algorithm is shown to scale linearly as well.
In~\prettyref{sec:discussion}, a parallelization study is performed and the performance in simulations evaluated.
Lastly, \prettyref{sec:conclusions} concludes the paper.

%%% Local Variables:
%%% mode: latex
%%% TeX-master: "block_smoothers"
%%% End:

%% file: block_smoothers_discretization.tex
\section{Spectral-element discretization}%

\subsection{The spectral-element method for hexahedral elements}%
\label{sec:discretization}
The \helmholtz\ equation in a domain~$\domain$ reads
\begin{align}
  \lambda \varu - \laplaceop \varu = \varf \label{eq:helmholtz_equation} \eqcomma
\end{align}
where~${\lambda}$ is a non-negative parameter, $\varu$ the solution, and~$\varf$ the right-hand side.
Discretizing the above equation leads to linear system of the form
\begin{align}
  \helmholtzop \vecu &= \vecF
\end{align}
with~$\helmholtzop$ as the discrete \helmholtz\ operator, $\vecu$ as solution vector and as $\vecF$ discrete right-hand side.
When using a spectral-element method with nodal basis functions, the global operator~$\helmholtzop$ consists of the assembly of element-local operators:
\begin{align}
  \helmholtzop &= \sum\limits_e\mathbf{Q}_e  \helmholtzop_{e} \mathbf{Q}_e^T \eqdot
\end{align}
In the above~${\helmholtzop_{e}}$ is the operator in element~$\domain_e$ and~$\mathbf{Q}_e$ is the respective assembly matrix~\cite{deville_2002_sem, lottes_2005_mg}.
For cuboidal tensor-product elements, the operators~$\helmholtzop_{e}$ simplify to
\begin{align}
  \begin{aligned}
    \helmholtzop_{e}
    &= \geomcoeff{0,e}\tp{\standardmassmat}{\standardmassmat}{\standardmassmat}\\
    &+ \geomcoeff{1,e}\tp{\standardmassmat}{\standardmassmat}{\standardstiffmat}\\
    &+ \geomcoeff{2,e}\tp{\standardmassmat}{\standardstiffmat}{\standardmassmat}\\
    &+ \geomcoeff{3,e}\tp{\standardstiffmat}{\standardmassmat}{\standardmassmat}
  \end{aligned}\label{eq:helmholtzop} \eqdot
\end{align}
Here, $\standardmassmat$ is the one-dimensional standard element mass matrix, approximated with GLL quadrature and~$\standardstiffmat$ the standard element stiffness matrix~\cite{deville_2002_sem}.
Both are of size~$(\poly+1)\times(\poly+1)$ where~$\poly$ is the polynomial degree.
Furthermore, the coefficients~$\geomcoeff{i,e}$ are given by
\begin{align}
  \mymat{\geomcoeff{}}_{e} &= \frac{h_{1,e}h_{2,e}h_{3,e}}{8}{\left(\lambda, \frac{4}{h_{1,e}^2}, \frac{4}{h_{2,e}^2}, \frac{4}{h_{3,e}^2} \right)}^{T}
\end{align}
where~${h_{i,e}}$ are the dimensions of element~$\Omega_e$ and $\otimes$ is the~{Kronecker} product operator, also referred to as tensor product~\cite{lynch_1964_tensors, deville_2002_sem}.
%
%Evaluating the element operator requires three one-dimensional matrix products per tensor-product operator, resulting in~$6{(\poly+1)}^4 + 2{(\poly+1)}^3$ operations when exploiting the diagonal mass matrix of the~GLL points.
%

\subsection{Static condensation}

The solution of~\eqref{eq:helmholtz_equation} depends on both the right-hand side and the boundary conditions, but not on values in the interior, which can be eliminated from the equation system.
This elimination is called static condensation, \schur\ complement or substructuring, and is typically employed to eschew interior degrees of freedom from the elements.
It generates an equation system with fewer degrees of freedom which is also better conditioned~\cite{couzy_1995_condensation}.

With static condensation the degrees of freedom are separated into element boundary ones, called~$\vecu_{\bound}$, and element interior ones, $\vecu_{\inner}$, such that~\eqref{eq:helmholtzop} can be rearranged to
\begin{align}
  \begin{pmatrix}
    \helmholtzop_{\bound\bound} & \helmholtzop_{\bound\inner}\\
    \helmholtzop_{\inner\bound} & \helmholtzop_{\inner\inner}
  \end{pmatrix}
                                  \begin{pmatrix}
                                    \vecu_{\bound}\\
                                    \vecu_{\inner}
                                  \end{pmatrix}
                                &=
                                  \begin{pmatrix}
                                    \vecF_{\bound}\\
                                    \vecF_{\inner}
                                  \end{pmatrix} \eqdot
  \intertext{The values on interior points equate to}
  \vecu_{\inner} &=  \helmholtzop_{\inner\inner}^{-1} \vecF_{\inner} - \helmholtzop_{\inner\inner}^{-1}\helmholtzop_{\inner\bound} \vecu_{\bound} \eqcomma
                   \intertext{leading to a system for the boundary degrees}
  \underbrace{\left(\helmholtzop_{\bound\bound} - \helmholtzop_{\bound\inner}\helmholtzop_{\inner\inner}^{-1}\helmholtzop_{\inner\bound}\right)}_{\condop}\underbrace{\vecu_{\bound}}_{\condu}
                                &= \underbrace{\vecF_{\bound}  - \helmholtzop_{\bound\inner}\helmholtzop_{\inner\inner}^{-1} \vecF_{\inner}}_{\condF}\\
  \Leftrightarrow      \condop \condu &= \condF \eqdot
                                        \intertext{When eliminating only element-interior degrees of freedom, the condensed operator consists of element operators, such that}
                                        \condop &= \sum\limits_e\mymat{\hat{Q}}_e  \condop_{e} \mymat{\hat{Q}}_e^T \eqdot
\end{align}
For cuboidal elements, the operators~$\condop_{e}$ can be written in tensor-product form and evaluated with just~$\order{\poly^3}$ operations~\cite{huismann_2017_condensation}.

\subsection{Additive~\schwarz\ methods}
\schwarz\ methods are a standard solution technique in continuous and discontinuous~\galerkin\ spectral-elements methods~\cite{feng_2001_mg, lottes_2005_mg, dryja_1994_substructuring, sherwin_2001_sem}.
Instead of solving the whole equation system, a~\schwarz\ method determines a correction of the current approximation by combining local results obtained from overlapping subdomains.
Repeating the process leads to convergence.

The current approximation of the solution, $\condu$, leaves a residual~$\condr$:
\begin{align}
  \condr &=  \condF - \condop \condu \eqdot
\end{align}
A correction~$\Delta\condu$ to~$\condu$ is sought so that the residual norm gets significantly lowered, i.e.\,${\condop \Delta\condu \approx \condr}$.
Computing the exact correction requires solution of~${\condop \condu = \condr}$, which is the overall target and expensive.
However, on small subdomains with~\dirichlet\ boundary conditions inversion is possible and relatively cheap.
Such local corrections~${\Delta \condu_{i}}$ are computed on multiple overlapping subdomains~$\Omega_{i}$
\begin{align}
  \condop_{i} \Delta \condu_{i} &= \condrestrict_{i}\condr
                                  \intertext{where~$\condrestrict_i$ is the Boolean restriction to subdomain~$\Omega_{i}$ and $\condop_{i}$ is the~\helmholtz\ operator restricted to the subdomain~$\Omega_{i}$}
                                  \condop_{i} &= \condrestrict_i \condop \condrestrict_i^T \eqdot
  \intertext{The correction is afterwards computed by suitably weighing the contributions from the subdomains:}
  \Delta \condu &= \sum_{i} \condrestrict^T_{i}\condweigh_{i} \underbrace{\condop_{i}^{-1}\condrestrict_{i} \condr}_{\Delta \condu_{i}} \label{eq:smoother_definition} \eqcomma
\end{align}
with~$\condweigh_{i}$ constituting the weight matrices for the subdomains.
Standard additive~\schwarz\ methods setting~$\condweigh_{i}$ to identity, only converge when using a relaxation factor~\cite{gander_2008_schwarz}.
Choosing~$\condweigh_{i}$ as the inverse multiplicity, i.e.\,the weighing each point with the number of subdomains he occurs in, lifts this restriction~\cite{lottes_2005_mg} and distance-based weight functions yield excellent convergence rates~\cite{stiller_2017_multigrid}.

With~\schwarz\ methods the choice of subdomain type and in particular the resulting overlap is paramount for the performance of the resulting algorithm~\cite{gander_2008_schwarz} and, hence, subdomains larger than one element are beneficial.
% The main component to computing the correction in the~\schwarz\ method~\eqref{eq:smoother_definition} is the inverse~\helmholtz\ operator on a subdomain.
%
In~\cite{haupt_2017_mg, haupt_2013_mg}, ${2^{d}}$ element blocks were utilized as subdomains in~$\mathbb{R}^d$, resulting in a vertex-based smoother.
Static condensation on this block leaves only the faces interconnecting the elements remaining, rather than all faces surrounding an element for an element-centered block.
This leads to considerably less degrees of freedom required to compute the inverse and renders the method referable for the condensed system compared to element-based subdomains.
As the block resembles a star in 2D it is called star smoother.

\begin{figure}[t]
  \hspace*{\fill}
  \includegraphics[width=0.25\textwidth]{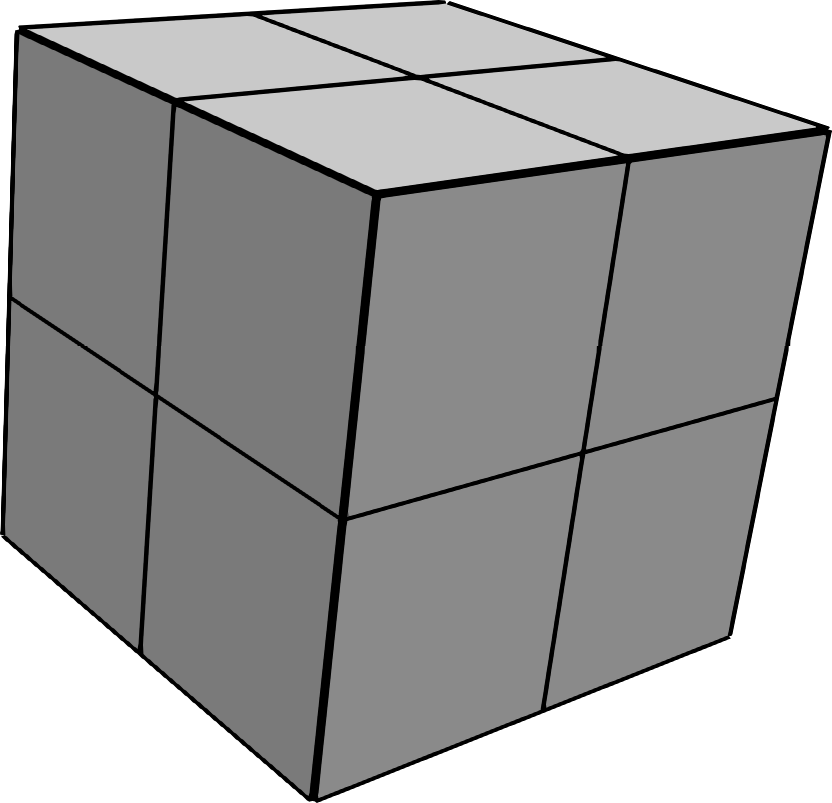}
  \hfill
  \includegraphics[width=0.25\textwidth]{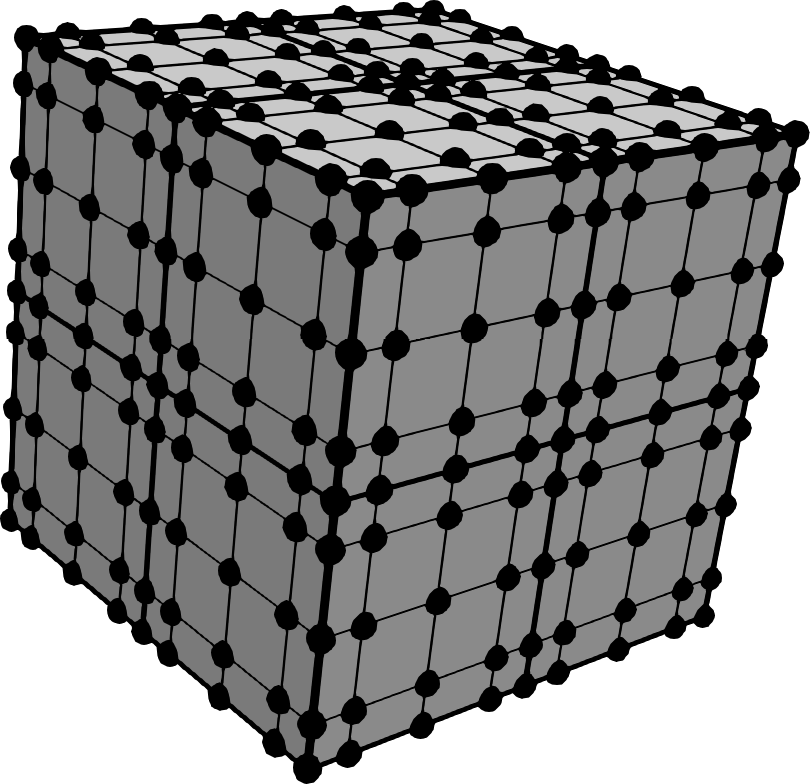}
  \hfill
  \includegraphics[width=0.25\textwidth]{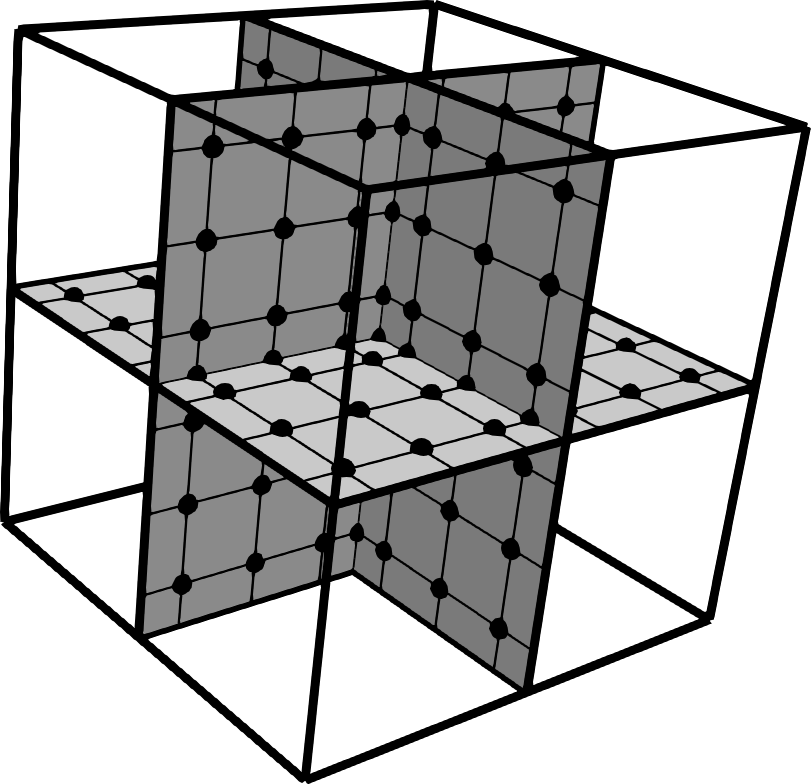}
  \hspace*{\fill}
  \caption{Block utilized for the star smoother. Left:~Elements of the subdomain, middle:~collocation points including~\textsc{Dirichlet} nodes, right:~collocation nodes in the condensed case.}%
  \label{fig:star}
\end{figure}

\prettyref{fig:star} depicts the~${2\times 2 \times 2}$ element block as well as the three-dimensional star.
As in a residual-based formulation the problem is homogeneous, the boundary nodes are dropped, resulting in~$\nstar = 2 \poly- 1$ data points per dimension.
For the condensed system the three planes connecting the elements remain, each with with~$\nstar^2$ degrees of freedom.
The resulting condensed star operator, $\condop_{i}$, can be decomposed into tensor-products, with the application utilizing~$\order{\nstar^3}$ operations when using the technique from~\cite{huismann_2017_condensation}.
However, as the operator consists of many suboperators, the inversion of it is not trivial when trying to retain the linear scaling.
For instance the implementation from~\cite{haupt_2017_mg} utilizes matrix inversion, leading to a dense matrix~$\condop_{i}^{-1}$, that couples every degree of freedom of the star with every other one, without capabilities for structure exploitation.
In three dimensions this results in~$2 \cdot {(3 \cdot \nstar^2)}^2 = 18 \cdot \nstar^4$ operations for application of the inverse.
%%% Local Variables:
%%% mode: latex
%%% TeX-master: "block_smoothers"
%%% End:

%% file: block_smoothers_smoothers.tex
\section{A linearly scaling block inverse}\label{sec:smoothers}

\subsection{Embedding the condensed system into the full system}

The~\schwarz\ method requires solution of the condensed system in subdomain~$\Omega_{i}$.
While the operator itself can be written using tensor products, finding an explicit inverse has so far eluded the community.
Matrix inversion, as proposed in~\cite{haupt_2017_mg} scales with~${\order{\nstar^4} = \order{\poly^4}}$ and, furthermore, is not matrix free by definition.
Overwhelming memory requirements are the result, which only allow for homogeneous meshes and limit the method to relatively low polynomial degrees.
These drawbacks need to be circumvented.
The main goal of this paper is to lift both stifling restrictions by deriving a matrix-free inverse as well as linearizing the operation count.
Achieving this naturally leads to a multigrid method scaling linearly with the number degrees of freedom, when increasing the number of elements as well as when increasing the polynomial degree.
To this end, the condensed system is embedded into the full equation system, where the fast diagonalization technique can be exploited to gain a matrix-free inverse, which will be factorized later on.

For embedding the stars into their~${2 \times 2 \times 2}$ element blocks considering one star is sufficient, hence, the subscript~$i$ is dropped in this discussion.
As done with the full equation system, the system is reordered into degrees on the star and remaining, element-internal degrees:
\begin{align}
  \begin{pmatrix}
    \helmholtzop_{\istar\istar} & \helmholtzop_{\istar\inner}\\
    \helmholtzop_{\inner\istar} & \helmholtzop_{\inner\inner}
  \end{pmatrix}
                                  \begin{pmatrix}
                                    \vecu_{\istar}\\
                                    \vecu_{\inner}
                                  \end{pmatrix}
                                &=
                                  \begin{pmatrix}
                                    \vecF_{\istar}\\
                                    \vecF_{\inner}
                                  \end{pmatrix}\\
  \Rightarrow \underbrace{\left(\helmholtzop_{\istar\istar} - \helmholtzop_{\istar\inner}\helmholtzop_{\inner\inner}^{-1}\helmholtzop_{\inner\istar}\right)}_{\condop}\vecu_{\istar}
                                &= \vecF_{\istar}  - \helmholtzop_{\istar\inner}\helmholtzop_{\inner\inner}^{-1} \vecF_{\inner} \eqdot
\end{align}
If the condensed system is embedded in the full system, the same solution is generated by the right-hand side of the condensed system.
The modified right-hand side
\begin{align}
  \starrhs &= \begin{pmatrix}\starrhs_{\istar}\\ \starrhs_{\inner}\end{pmatrix} = \begin{pmatrix}\vecF_{\istar} - \helmholtzop_{\istar\inner}\helmholtzop_{\inner\inner}^{-1}\vecF_{\inner}\\ 0\end{pmatrix}\label{eq:star_right_hand_side_full}
  \intertext{infers a solution~$\starvar$ satisfying}
  \begin{pmatrix}
    \helmholtzop_{\istar\istar} & \helmholtzop_{\istar\inner}\\
    \helmholtzop_{\inner\istar} & \helmholtzop_{\inner\inner}
  \end{pmatrix}
                                  \begin{pmatrix}
                                    \starvar_{\istar}\\
                                    \starvar_{\inner}
                                  \end{pmatrix}
           &=
             \begin{pmatrix}
               \starrhs_{\istar}\\
               0
             \end{pmatrix}\\
  \Rightarrow \underbrace{\left(\helmholtzop_{\istar\istar} - \helmholtzop_{\istar\inner}\helmholtzop_{\inner\inner}^{-1}\helmholtzop_{\inner\istar}\right)}_{\condop}\starvar_{\istar}
           &= \starrhs_{\istar} = \vecF_{\istar}  - \helmholtzop_{\istar\inner}\helmholtzop_{\inner\inner}^{-1} \vecF_{\inner} \eqdot
\end{align}
The operator~${\condop}$ is positive definite and, hence, the linear system possesses a unique solution, i.e.~$\starvar_{\istar} = \vecu_{\istar}$.
However, the solution of the whole operator is unique as well.
As the right-hand sides differ, the solutions differ in the interior, i.e.~$\starvar_{\inner} \ne \vecu_{\inner}$, which is not required for the condensed case.

Solving using the full operator with~$\starrhs$ leads to the correct solution on the star, enabling the usage of solution methods from the full system to attain the solution in the condensed one.

\subsection{Tailoring fast diagonalization for static condensation}
As the full system replicates the solution of the condensed system, a method to attain the solution for the condensed system can be obtained by restriction from the full system.
This generates the necessity to investigate the operator of the full star.
As with the element~\helmholtz\ operator, using lexicographic ordering leads to the collocation nodes on the full star~$\domain_{i}$ exhibiting a tensor-product structure, allowing the operator~$\helmholtzop_{i}$ to be formulated as in~\eqref{eq:helmholtzop}:
\begin{align}
  \begin{aligned}
    \helmholtzop_{i}
    &= \geomcoeff{0}\tp{\massmat_{i}}{\massmat_{i}}{\massmat_{i}}\\
    &+ \geomcoeff{1}\tp{\massmat_{i}}{\massmat_{i}}{\stiffmat_{i}}\\
    &+ \geomcoeff{2}\tp{\massmat_{i}}{\stiffmat_{i}}{\massmat_{i}}\\
    &+ \geomcoeff{3}\tp{\stiffmat_{i}}{\massmat_{i}}{\massmat_{i}} \eqdot
  \end{aligned}\label{eq:helmholtzopstar}
\end{align}
Here~$\massmat_{i}$ and~$\stiffmat_{i}$ are the one-dimensional mass and stiffness matrices restricted to full star~$\domain_{i}$.
They correspond to the one-dimensional matrices assembled from the elements, as shown in~\prettyref{fig:star}.
For sake of readability, the same stiffness and mass matrices were utilized in all three dimensions.
When varying element dimensions inside a block, these matrices will differ.

The inverse of the block~\helmholtz\ operator can be expressed via the fast diagonalization~\cite{lynch_1964_tensors}.
A generalized eigenvalue decomposition of~${\transmat_{i}}$ with respect to~${\massmat_{i}}$ is performed, such that
\begin{subequations}
\begin{align}
  \transmat_{i}^{T} \stiffmat_{i} \transmat_{i} &= \eigenvaluemat_{i}\\
  \transmat_{i}^{T} \massmat_{i} \transmat_{i} &= \identitymat \eqdot
\end{align}
\end{subequations}
The matrix~$\transmat_{i}$ is a non-orthogonal transformation matrix and~$\eigenvaluemat_{i}$ is the diagonal matrix comprising the generalized eigenvalues of~$\stiffmat_{i}$.
This decomposition allows to express the inverse of~\eqref{eq:helmholtzopstar} as
\begin{align}
  \helmholtzop_{i}^{-1} &= \ptp{\transmat_{i}}{\transmat_{i}}{\transmat_{i}} \diagonal_{i}^{-1} \ptp{\transmat^{T}_{i}}{\transmat_{i}^{T}}{\transmat_{i}^{T}} \label{eq:fast_diagonalization} \eqcomma
                          \intertext{where the three-dimensional eigenvalues are stored in the diagonal matrix}
                          \diagonal_{i} &= \geomcoeff{0} \ptp{\identitymat}{\identitymat}{\identitymat}
                                          + \geomcoeff{1}\ptp{\identitymat}{\identitymat}{\eigenvaluemat_{i}}
                                          + \geomcoeff{2}\ptp{\identitymat}{\eigenvaluemat_{i}}{\identitymat}
                                          + \geomcoeff{3}\ptp{\eigenvaluemat_{i}}{\identitymat}{\identitymat} \eqdot
\end{align}
The application of the above operator~\eqref{eq:fast_diagonalization} with tensor-products utilizes~$12\cdot\nstar^4$ operations, which still scales super-linearly with the number of degrees of freedom when increasing~$\poly$.

To attain the solution of the condensed star, using the reduced right-hand side from~\eqref{eq:star_right_hand_side_full} suffices.
Applying~\eqref{eq:fast_diagonalization} consists of three steps: Mapping into the three-dimensional eigenspace, applying the inverse eigenvalues, and then mapping back.
Due to~$\starrhs_{i}$ being zero in element-interior regions, the three-dimensional tensor-product operators in~\eqref{eq:fast_diagonalization} only operate on the three planes of the star.
Hence, computing~${\starrhs_{\eigenspace} = \ptp{\transmat^{T}_{i}}{\transmat_{i}^{T}}{\transmat_{i}^{T}} \starrhs_{i}}$ can be decomposed into the application of three separate tensor-products working on two-dimensional data instead of one working on the whole data set.
These two-dimensional tensor products require~$\order{\nstar^3}$ operations when expanding to the star eigenspace last.
Thus, the operation of mapping the right-hand side~$\starrhs$ into the star eigenspace scales linearly when increasing~$\poly$.
Applying the diagonal matrix of inverse eigenvalues is linear as well.
Lastly, the result is only required on faces of the star, allowing not to compute interior degrees, resulting in~$\order{\nstar^3}$ operations for mapping from the eigenspace to the star planes.
The combination of all three operations yields an inverse that can be applied with linear complexity.

\prettyref{alg:inverted_star_tensor} summarizes the resulting procedure to apply the inverse full star in the condensed system.
The right-hand side is extracted and stored on the three faces of the star.
These are perpendicular to the~$x_1$, $x_2$, and~$x_3$ directions and called~$\face_{1}$, $\face_{2}$, and~$\face_{3}$, respectively.
The data on each face is stored as matrix of extent~$\nstar \times \nstar$, where the indices ranging from~$\inner = \{-\poly \dots \poly\}$ with index~$0$ corresponding to the 1D index of the face in the full system.
First, the inverse multiplicity~$\bm{\mathcal{M}}^{-1}$ is applied to account for multiply stored values.
Then the two-dimensional transformation matrix~$\tptwo{\transmat^{T}}{\transmat^{T}}$ transforms the three faces.
The results are mapped into the star eigenspace using permutations of~$\tp{\identitymat}{\identitymat}{\transmat^{T}_{\inner 0}}$.
There, the inverse eigenvalues are applied and, lastly, the reverse order of operations is utilized to map back to the star faces.

Overall, \prettyref{alg:inverted_star_tensor} requires~$18$ one-dimensional matrix products, 9 to expand from two dimensions to three and further~9 to reduce from three-dimensional to two-dimensional data.
Moreover, one multiplication in eigenspace is required, leading to~${37 \nstar^3}$ operations required to apply the algorithm.
\begin{algorithm}
  \caption{Calculation of the solution~$\starvar$ on the three faces of the star, $\face_{1}$, $\face_{2}$, and $\face_{3}$, which are perpendicular to~$x_1$, $x_2$, and~$x_3$, respectively, from the star residual~$\starrhs$ stored on these faces.
    Values on multiply occurring data points in the residual are divided by their multiplicity, allowing to retain the tensor-product structure of the operator.}\label{alg:inverted_star_tensor}
  \begin{algorithmic}
    \State $\starrhs_{\face} \gets\bm{\mathcal{M}}^{-1}\starrhs_{\face}$ \Comment{account for multiply occurring degrees of freedom}
    \State $\starrhs_{\eigenspace} \gets\phantom{+} \ptp{\transmat^{T}}{\transmat^{T}}{\transmat_{\inner 0}^{T}} \starrhs_{\face_1}$ \Comment{contribution from face perpendicular to~$x_1$}
    \State $\phantom{\starrhs_{\eigenspace} \gets}  +     \ptp{\transmat^{T}}{\transmat_{\inner 0}^{T}}{\transmat^{T}} \starrhs_{\face_2}$\Comment{contribution from face perpendicular to~$x_2$}
    \State $\phantom{\starrhs_{\eigenspace} \gets}  +     \ptp{\transmat^{T}_{ \inner 0}}{\transmat^{T}}{\transmat^{T}} \starrhs_{\face_3}$\Comment{contribution from face perpendicular to~$x_3$}
    \State $\starvar_{\eigenspace} \gets \diagonal^{-1} \starrhs_{\eigenspace}$ \Comment{application of inverse in the eigenspace}
    \State $\starvar_{\face_{1}} \gets \ptp{\transmat}{\transmat}{\transmat_{0 \inner}} \starvar_{\eigenspace}$ \Comment{map solution from eigenspace to~$x_1$ face}
    \State $\starvar_{\face_{2}} \gets \ptp{\transmat}{\transmat_{0 \inner}}{\transmat} \starvar_{\eigenspace}$ \Comment{map solution from eigenspace to~$x_2$ face}
    \State $\starvar_{\face_{3}} \gets \ptp{\transmat_{0 \inner}}{\transmat}{\transmat} \starvar_{\eigenspace}$ \Comment{map solution from eigenspace to~$x_3$ face}
  \end{algorithmic}
\end{algorithm}

\begin{figure}
  \centering
  \includegraphics{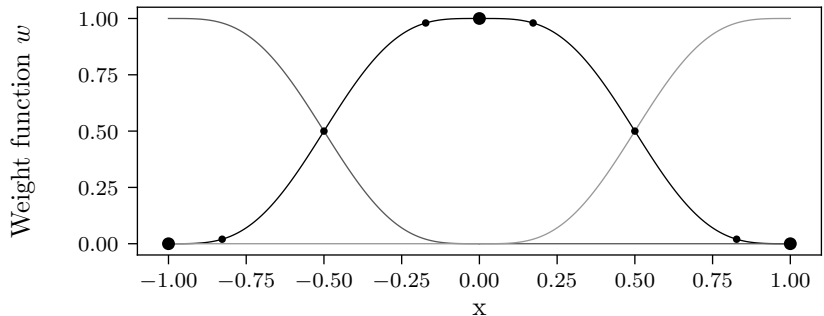}
  \caption{%
    Polynomial weight functions of degree~${\poly_{\weigh} = 7}$ on two adjoining elements in one dimension.
    The polynomials are one on their respective vertex and zero on every other one, leading to a partition of one.
  }%
  \label{fig:weights}
\end{figure}

\prettyref{alg:smoother} shows the resulting~\schwarz-type smoother.
For every vertex, the data of the surrounding elements is gathered.
The inverse is applied, leading to a correction~$\Delta \condu_{i}$ on the star.
Afterwards, the contribution of each vertex is weighed to gain the correction in the surrounding elements, where~$\condweigh_{i}$ is the weight matrix for each star.
It is constructed by restricting the tensor-product of one-dimensional diagonal weight matrices~$\weigh = \tp{\weigh^{\mathbf{1D}}}{\weigh^{\mathbf{1D}}}{\weigh^{\mathbf{1D}}}$ to the condensed system.
As in~\cite{stiller_2017_multigrid}, these are populated by smooth polynomials of degree~$\poly_{\weigh}$, which are one on the vertex of the star and zero on all other vertices.
For polynomial degrees larger than one, higher derivatives are zero at the vertices smoothing the transition as shown in~\prettyref{fig:weights}.
In the studies shown below~$\poly_{\weigh}=7$ was utilized.
\begin{algorithm}
  \caption{Smoother smoothing over the stars corresponding to the~$n_{\mathrm{V}}$ vertices.}%
  \label{alg:smoother}
  \begin{algorithmic}
    \Function{Smoother}{$\condr$}
    \For{$i = 1, n_{\mathrm{V}}$}
    \State{$\condrhs_{i} \gets  \bm{\mathcal{\hat{M}}}_{i}^{-1} \condrestrict_{i}\condr $}
    \Comment{extraction of data}
    \State{$\Delta \condvar_{i} \gets \condop^{-1}_{i}\condrhs_{i}$}
    \Comment{inverse on stars}
    \EndFor{}
    \State{$\Delta \condu \gets \sum_{i=1}^{n_{\mathrm{V}}}\condweigh_{i}\condrestrict_{i}^{T} \Delta \condu_{i}$}
    \Comment{star contributions from~$8$ vertices per element}
    \Return{$\Delta \condu$}
    \EndFunction{}
  \end{algorithmic}
\end{algorithm}

\subsection{Implementation of boundary conditions}
The additive~\schwarz\ method works on subdomains corresponding to the eight elements surrounding a vertex.
This poses a problem for the implementation on boundary vertices:
The subdomain reduces in size, depending on the number of boundary conditions present.
In the worst case~${ 5^3 = 125 }$ different implementations of the inverse are required when accounting for~\dirichlet\ and \neumann\ boundary conditions.
This section seeks an approach using the same treatment for boundary vertices as for interior ones, lowering the implementation effort.
With structured grids and parallelization, one layer of ghost elements is assumed.
The approach changes the matrices for the stars, rather than their structure.
Due to the structured mesh, investigating the one-dimensional case suffices, as the operator retains its tensor-product structure.

\begin{figure}
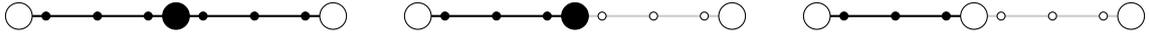

  \includegraphics[width = 0.3 \textwidth]{images/boundary_condition_full.tikz}
  \hfill
  \includegraphics[width = 0.3 \textwidth]{images/boundary_condition_neumann.tikz}
  \hfill
  \includegraphics[width = 0.3 \textwidth]{images/boundary_condition_dirichlet.tikz}
  \caption{Implementation of boundary conditions for the one-dimensional case on the right domain boundary.
    Utilized data points are drawn in black, non-utilized ones in white.
    Left: Used smoother block consisting of two elements of polynomial degree~${\poly= 4}$.
    Middle: A \neumann\ boundary condition decouples the right element.
    Right: Homogeneous \dirichlet\ boundary condition decouples the middle vertex as well.}%
  \label{fig:boundary_conditions}
\end{figure}

\prettyref{fig:boundary_conditions} depicts the boundary condition implementation in one dimension.
In the periodic case the ghost elements stay part of the domain, requiring no change in the operator.
However, with~\neumann\ boundaries, the ghost element is not part of the domain anymore and needs to be decoupled.
In effect, the one-dimensional domain around the vertex shrinks to one element.
To account for this, the corresponding transformation matrices and eigenvalues are padded with identity, so that the same operator size is used throughout the whole domain.
Furthermore the discrete right-hand side is zeroed outside of the domain.
A similar treatment can be used for~\dirichlet\ boundaries.
In addition to decoupling the ghost element, the boundary vertex is decoupled and the corresponding right-hand side zeroed.
As the decoupled points only map onto themselves, the computed correction is zero, leading to no correction for~\dirichlet\ boundary conditions which keep their former value.

For structured grids in multiple dimensions, the operator consists of a tensor-product of the matrices for the one-dimensional case.
Four cases exist:
In the first case every element is part of the domain, requiring no change.
In the second one, a boundary condition is present in one dimension, changing one set of matrices.
In the third case, boundary conditions are present in two dimensions, requiring change in two sets of matrices.
In the fourth case, boundary conditions are present in all directions, leading to an equation system for the point alone.
Padding the respective transformation matrices with identity leads to the decoupled parts only mapping onto themselves.
%
%Further tensor products only work in the corresponding half-planes.
%
Hence, no information is transported from outside the domain into it, generating a correct implementation of boundary conditions.

\subsection{Extension to element-centered block smoothers}\label{sec:element_centered}
So far only vertex-based smoothers were considered, which only allow for half an element of overlap with the~$2^3$ element cube.
When increasing the overlap to one element, element-centered subdomains are more favorable.
Furthermore many algorithms in the literature work on element-centered subdomains~\cite{lottes_2005_mg, stiller_2017_multigrid} and having a drop-in replacement is beneficial.

\begin{figure}[t]
  \hspace*{\fill}
  \includegraphics[width=0.4\textwidth]{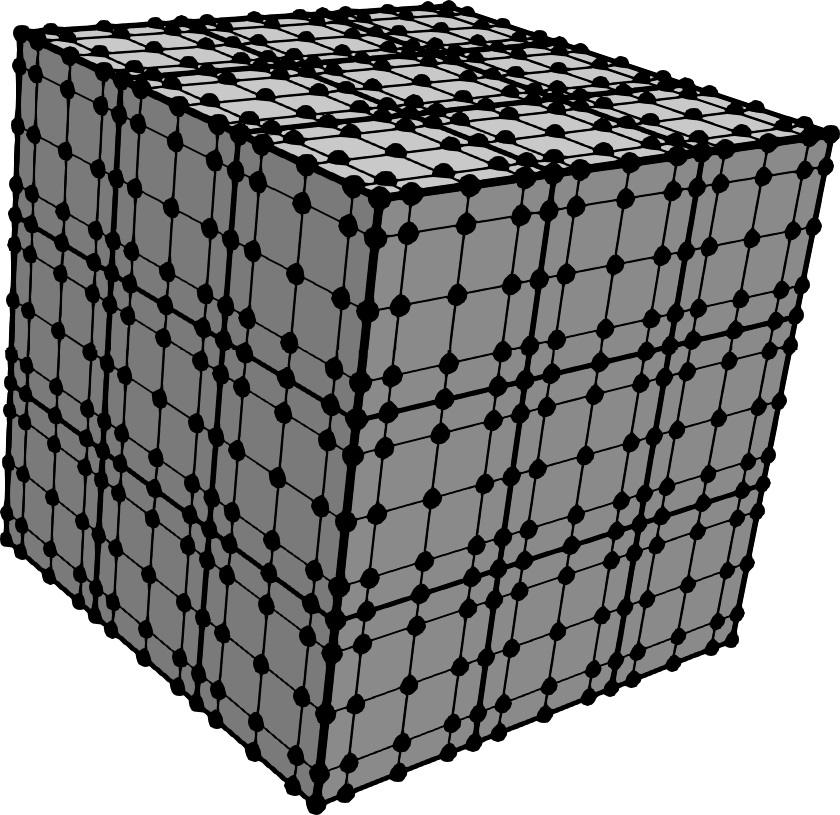}%
  \hfill
  \includegraphics[width=0.4\textwidth]{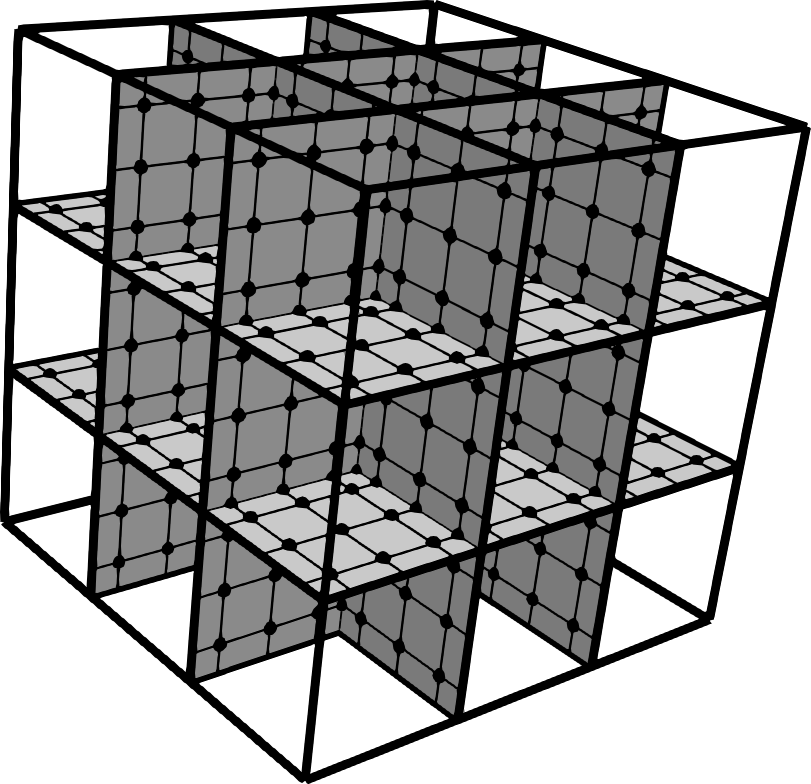}%
  \hspace*{\fill}\\
  \caption{Block utilized for the element-centered smoothers.
    Left:~full system including \textsc{Dirichlet} nodes.
    Right:~Collocation nodes corresponding to the condensed setup.}%
  \label{fig:element_centered_block}
\end{figure}

\prettyref{fig:element_centered_block} depicts an element-centered subdomain overlapping into the neighboring elements.
Compared to the star smoother, three elements per direction are required.
While the amount of overlap can be chosen depending on the polynomial degree, the simplest case is the full overlap, which is utilized here for demonstration.
As with the vertex-centered case, the operator in the block can be written in tensor-product form~\eqref{eq:helmholtzop}, with the one-dimensional matrices being replaced by those restricted to the subdomain.
Again fast diagonalization is applicable.

With the techniques used for condensing and embedding the star into the full system, the element-centered block can be condensed and embedded.
Hence, solution techniques from the full~$3^d$ element block can be restricted to the condensed element-centered block shown in~\prettyref{fig:element_centered_block}.
Furthermore, the same arguments can be applied for attaining a linearly-scaling inverse from the fast diagonalization.
Due to six faces being present instead of three, the number of operations increases to~${73 \nstar^3}$, where with the full overlap~${\nstar = 3 \poly - 1}$.
However this needs to compare to the fast diagonalization in the full system, which requires~${12 \nstar^4}$, rendering the new algorithm more efficient starting from~${\poly = 4}$.

%%% Local Variables:
%%% mode: latex
%%% TeX-master: "block_smoothers"
%%% End:

%% file: block_smoothers_multigrid.tex
\section{Multigrid method}
\subsection{Multigrid Algorithm}

The $\poly$-multigrid is a well-researched building block for higher-order methods~\cite{ronquist_1987_mg, lottes_2005_mg, pasquetti_2009_multigrid}.
Where traditional~$h$-multigrid coarsens the mesh width~$h$ consecutively~\cite{brandt_1982_mg, bramble_1995_mg, hackbusch_1985_mg}, $\poly$-multigrid lowers the polynomial degree for coarser grids.
Here, the polynomial degree is reduced from level~$L$ to level~$0$, with the series of polynomial degrees defined as
\begin{subequations}
  \begin{align}
    \forall\ 0 \leq  l < L:\quad \poly_{l} &= \poly_0 \cdot 2^{l}\\
    \poly_{L} &= \poly \eqdot
  \end{align}
\end{subequations}

The implementation requires four main components:
An operator calculating the residual, a smoother smoothing out the high-frequency components on the current level, grid transfer operators, and a coarse grid solver for the equation system at~$\poly_0$.
In the present case, the system operator is the condensed operator in tensor-product form~\cite{huismann_2017_condensation}.
The prolongation from level~$l-1$ to level~$l$, ${\condinterpol_{l}}$, is the embedded interpolation, restricted from the tensor-product of one-dimensional prolongation operators to the condensed system, and the restriction operator is its transpose.
Lastly, the coarse grid solver considered in this paper is a conjugate gradient solver working in the condensed system at~${\poly_{0} = 2}$.

A standard V-cycle, as shown in~\prettyref{alg:vcycle}, is used.
It allows for a varying number of pre- and post-smoothing steps~$\nu_{\mathrm{pre},l}$ and~$\nu_{\mathrm{post},l}$ on level~$l$.
Two methods of choosing these are considered here.
In the first one, one pre- and one post-smoothing step are used, the second one increases both with each level~${\nu_{\mathrm{pre},l} = \nu_{\mathrm{post},l} = 2^{L-l}}$, stabilizing the method for non-uniform meshes~\cite{stiller_2017a_multigrid}.

\begin{algorithm}
  \caption{Multigrid V-cycle for the condensed system using~$\nu_{\mathrm{pre}}$ pre- and~$\nu_{\mathrm{post}}$ post-smoothing step}%
  \label{alg:vcycle}
  \begin{algorithmic}
    \Function{MultigridCycle}{$\condu$, $\condF$}
    \State{$\condu_{L} \gets \condu$}
    \State{$\condF_{L} \gets \condF$}
    \For{$l = L, 1, -1$}
    \If{ $l \neq L$}
    \State $\condu_l \gets 0$
    \EndIf
    \For{$ i = 1, \nu_{\mathrm{pre},l}$}
    \State $\condu_l \gets \condu_{l} + \textsc{Smoother}( \condF_l - \condop_{l} \condu_{l})$    \Comment{Presmoothing}
    \EndFor

    \State $\condF_{l-1} \gets \condinterpol_{l-1}^{T} \left(\condF_l - \condop_l \condu_l \right)$
    \Comment{Restriction of residual}
    \EndFor{}
    \State $\mathrm{Solve}(\condop_{0} \condu_{0} = \condF_0)$
    \Comment{Coarse grid solve}
    \For{$l = 1, L$}
    \State $\condu_l \gets \condu_l + \condinterpol_{l}\condu_{l-1}$
    \Comment{Prolongation of correction}
    \For{$ i = 1, \nu_{\mathrm{post},l}$}
    \State $\condu_l \gets \condu_{l} + \textsc{Smoother}( \condF_l - \condop_{l} \condu_{l})$    \Comment{Postsmoothing}
    \EndFor
    \EndFor

    \Return $\condu \gets \condu_{L}$
    \EndFunction
  \end{algorithmic}
\end{algorithm}

\prettyref{alg:mg} shows the resulting multigrid algorithm for the condensed system.
First, the variables from the full system are condensed: Solution and right-hand side are restricted, with the latter gaining contributions from the inner element.
V-cycles are  performed until convergence is attained.
Afterwards, the solution in the inner elements is computed.
\begin{algorithm}
  \caption{Multigrid algorithm for the condensed system}%
  \label{alg:mg}
  \begin{algorithmic}
    \Function{MultigridSolver}{$\vecu = {(\vecu_{\bound}, \vecu_{\inner})}^T$, $\vecF = {(\vecF_{\bound}, \vecF_{\inner})}^T$}
    \State{$\condu \gets \vecu_{\bound}$}
    \Comment{restrict to element boundaries}
    \State{$\condF \gets \vecF_{\bound} - \helmholtzop_{\bound\inner}\helmholtzop_{\inner\inner}^{-1}\vecF_{\inner}$}
    \Comment{restrict and condense inner element RHS}
    \While{$\sqrt{\condr^T \condr} > \varepsilon$}

    \State{$\condu \gets \textsc{MultigridCycle}(\condu, \condF)$}
    \Comment{fixpoint iteration with multigrid cycle}
    \EndWhile{}
    \State{$\vecu \gets {(\condu, \helmholtzop_{\inner\inner}^{-1}(\vecF_{\inner} - \helmholtzop_{\inner\bound}\condu ) )}^T $}
    \Comment{regain interior degrees of freedom}
    \Return{$\vecu$}
    \EndFunction{}
  \end{algorithmic}
\end{algorithm}

To enhance multigrid performance with inhomogeneous meshes, we consider so-called \krylov\ acceleration~\cite{oosterlee_1998_mg, woodward_1998_newton}.
Instead of as a direct solution method, the multigrid cycle serves as preconditioner in a pCG method.
However, due to the weighing, the smoother is not symmetric and, hence, standard~pCG not guaranteed to converge~\cite{hestenes_1952_cg}.
Following~\cite{stiller_2016_multigrid}, the inexact preconditioned~CG method is utilized as remedy for non-symmetric preconditioners~\cite{golub_1999_ipcg}.
\prettyref{alg:cgmg} shows the \textsc{Krylov}-accelerated multigrid algorithm for the condensed system.
\begin{algorithm}
  \caption{Krylov-accelerated multigrid algorithm for the condensed system}%
  \label{alg:cgmg}
  \begin{algorithmic}
    \Function{ipCGMultigridSolver}{$\vecu = {(\vecu_{\bound}, \vecu_{\inner})}^T$, $\vecF = {(\vecF_{\bound}, \vecF_{\inner})}^T$}
    \State{$\condu \gets \vecu_{\bound}$}
    \Comment{restrict to element boundaries}
    \State{$\condF \gets \vecF_{\bound} - \helmholtzop_{\bound\inner}\helmholtzop_{\inner\inner}^{-1}\vecF_{\inner}$}
    \Comment{restrict and condense inner element RHS}
    \State{$\condr \gets \condF - \condop \condu$}
    \Comment{initial residual}
    \State{$\conds \gets \condr$} \Comment{ensures $\beta = 0$ on first iteration}
    \State{$\condp \gets 0$}\Comment{initialization}
    \State{$\delta = 1$}\Comment{initialization}
    \While{$\sqrt{\condr^T \condr} > \varepsilon$}
    \State{$\condz \gets \textsc{MultigridCycle}(0, \condr)$}
    \Comment{preconditioner}
    \State{$\gamma \gets \condz^T \condr$}
    \State{$\gamma_0 \gets \condz^T \conds$}
    \State{$\beta = (\gamma - \gamma_0) / \delta$}
    \State{$\delta = \gamma$}
    \State{$\condp \gets \beta \condp +  \condz$} \Comment{update search vector}
    \State{$\condq \gets \condop \condp$} \Comment{compute effect of~$\condp$}
    \State{$\alpha = \gamma / (\condq ^T \condp)$} \Comment{compute step width}
    \State{$\conds \gets \condr$} \Comment{save old residual}

    \State{$\condu \gets \condu + \alpha \condp$} \Comment{update solution}
    \State{$\condr \gets \condr - \alpha \condq$} \Comment{update residual}
    \EndWhile{}
    \State{$\vecu \gets {(\condu, \helmholtzop_{\inner\inner}^{-1}(\vecF_{\inner} - \helmholtzop_{\inner\bound}\condu ) )}^T $}
    \Comment{Regain interior degrees of freedom}
    \Return{$\vecu$}
    \EndFunction:
  \end{algorithmic}
\end{algorithm}

\subsection{Complexity of resulting algorithms}
Both solvers presented in the previous section consist of three phases:
First, the condensation of the right-hand side, requiring~$W_{\mathrm{pre}}$ operations, then, the solution process with~$W_{\mathrm{sol}}$ operations, and, lastly, the recovery of interior degrees of freedom with~$W_{\mathrm{post}}$ operations.

Both, condensation and recomputation of inner degrees of freedom involve mapping from the current three-dimensional basis into the~\helmholtz\ operator eigenspace of the element and vice-versa.
%
% When implemented with a tensor product, both scale with~$\order{\poly^4\nelement}$.
% %
% The above can be reduced to~$\order{\poly^3}$ when using a basis comprising the eigenfunctions of the static condensation, however this is only prudent for linear problems, as, e.g., convection terms need to be computed on a nodal basis.
% %
% Hence, the coordinate transformation would only move to a different part of the code.
%
Hence, the total work required for pre- and post-processing is~$W_{\mathrm{pre}} + W_{\mathrm{post}} = \order{\poly^4 \nelement}$.
The iteration process consists, mainly, of applying the multigrid cycle.
The smoother on the fine grid scales with~$\order{\poly^3}$ as well as the operator~\cite{huismann_2017_condensation}.
For coarser grids and a constant number of smoothing steps, the cost of operator and smoothing scales with~${(1 / 2)}^3 = 1/8$ per coarsening, leading to a geometric series limited by a factor of~$8/7$.
Using~${\nu_{\mathrm{pre},l} = \nu_{\mathrm{post},l} = 2^{L-l}}$ leads to a lowering of the effort by~${2 \cdot (1 / 2)}^3 = 1/4$ and, hence, a limit of~$4/3$.
Lastly, the cost of the coarse grid solver scales with~$\order{p_0^3 \nelement^{\alpha}} = \order{\nelement^{\alpha}}$, where~$\alpha$ depends on the solution method.
For the examples presented below a method with~${\alpha=4/3}$ is employed.
We remark, however, that reaching~$\alpha=1$ can be attained by reverting to an appropriate low-order multigrid solver~\cite{brandt_1982_mg}.
Hence, the overall cost per cycle is~$\order{\poly^3\nelement}$ and the work required for a multigrid solution of the~\helmholtz\ equation is
\begin{align}
  W_{\mathrm{total}} &= W_{\mathrm{pre}} + N_{\mathrm{cycle}}W_{\mathrm{cycle}} + W_{\mathrm{post}}\\
  \Rightarrow W_{\mathrm{total}} &= \order{\poly^4\nelement} + N_{\mathrm{cycle}}\order{\poly^3 \nelement}
\end{align}
In the above, the main contributions stem form the fine grid costs:
Two smoothing steps with~$2 \cdot 37 \nstar ^3 \approx 2 \cdot 37 (2\poly^3) = 592 \poly^3$ operations per cycle, whereas pre- and post-processing require~$12 \poly^4 + 25 \poly^3$ each.
When assuming two multigrid cycles, the cost of the smoother dominates the cost of pre- and post-processing until~${\poly > 48}$.
Hence, for all practical purposes, the multigrid algorithm scales linearly with the number of elements, as well as with the polynomial degree.

%%% Local Variables:
%%% mode: latex
%%% TeX-master: "block_smoothers"
%%% End:

%% file: block_smoothers_operator_runtimes.tex
\section{Results}\label{sec:results}

\subsection{Runtimes for the  star inverse}
To test the efficiency of the inversion operator on the star, runtime tests were conducted.
Three variants were considered: A tensor-product variant using fast diagonalization in the full system called~``TPF'', a version using matrix-matrix multiplication in the condensed system to apply the  precomputed inverse called~``MMC'', and the proposed tensor-product variant implementing the fast diagonalization in the condensed system, called~``TPC''.
These were implemented in Fortran~2008 using double precision.
The tensor-product variants use loops, with the outermost one corresponding to the subdomain, leading to inherent cache blocking.
As~${\nstar = 2 \poly - 1}$ is odd by definition, the operator was padded by one, inferring an even operator size which, in combination with compile-time known loop bounds, enables compiler optimization.
Furthermore, the single-instruction multiple-data~(SIMD) compiler directive~\texttt{!dir\$ simd}, which is specific to the~Intel compiler, enforced vectorization of the inner-most non-reduction loop.
The variant ``MMC'' consisted of one call to DGEMM from BLAS, with the inverse being only computed for one subdomain and utilized for all.

A single node of the high-performance computer Taurus at ZIH Dresden served as test platform.
It comprised two~Xeon E5-2590 v3 processors, with twelve cores each, running at~${2.5\ \mathrm{GHz}}$.
Only one core was utilized, allowing for~${40\ \mathrm{GFLOP/s}}$ as maximum performance~\cite{hackenberg_2015_energy}.
The Intel Fortran Compiler v.\,2018 compiled the programs with the corresponding Intel Math Kernel Library~(MKL) serving as BLAS implementation.

\begin{figure}[h]
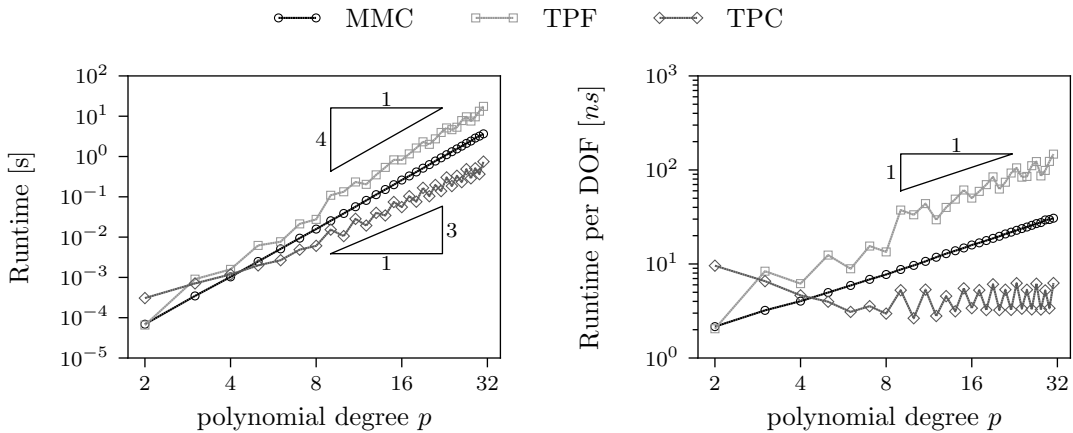

  \hspace*{\fill}
  \includegraphics{build/plot_runtime_test_block_inverse_legend.pgf}
  \hspace*{\fill}\\
  \hspace*{\fill}
  \includegraphics{build/plot_runtime_test_block_inverse_times.pgf}
  \hfill
  \includegraphics{build/plot_runtime_test_block_inverse_times_per_dof.pgf}
  % \hspace*{\fill}
  % \includegraphics{build/plot_runtime_test_block_inverse_gflops.pgf}
  \caption{Results for the application of the inverse star operator.
    Left: Operator runtimes when using the fast diagonalization in the full system~(TPF), applying the inverse via a matrix-matrix product in the condensed system~(MMC), and using the inverse via tensor-product factorization~(TPC).
    Right:~Runtimes per equivalent number of degrees of freedom~(DOF) being computed as~${(2\poly-1)}^3$ per block.}%
  \label{fig:runtime_measurements}
\end{figure}

The operators were used on~$500$ star subdomains, corresponding to 500 vertices being present on the process.
For these, the polynomial degree varied from~${\poly=2}$ to~${\poly=32}$.
Application of the operators was repeated 101 times, with the time last 100 times being measured via~MPI\_Wtime to preclude instantiation effects.
\prettyref{fig:runtime_measurements} depicts the runtime of the three variants.
All three exhibit the expected slopes: The fast diagonalization in the full system and the matrix-matrix variant have slope four, whereas the fast diagonalization in the condensed system has slope three.
This leads to the runtime per degree of freedom of the former two growing, whereas the proposed variant attains a constant runtime per degree of freedom for~${\poly > 8}$, with the only difference lying between even and odd polynomial degrees.
For even~$\poly$ the loop size is a multiple of the SIMD size of~$4$ leading to~${12\ \mathrm{GFLOP/s}}$, whereas for odd~$\poly$ the performance degrades by a factor of 2.
While for low~$\poly$ the matrix-matrix variant is the fastest, with the tensor-product variant in the full system being slightly slower, the fast diagonalization variant in the condensed system is faster starting from~$\poly = 5$, outpacing ``TPF'' by two orders of magnitude for~$\poly=32$.
%
%%% Local Variables:
%%% mode: latex
%%% TeX-master: "block_smoothers"
%%% End:

%% file: block_smoothers_solver_tests.tex
\subsection{Solver runtimes  for homogeneous meshes}%
\label{sec:solver_runtimes}
To verify that the new multigrid solvers scale linearly with the number of degrees of freedom, runtime tests were conducted utilizing the testcase from~\cite{huismann_2017_condensation}.
The \helmholtz\ problem is considered in a domain~${\domain = {(0,2\pi)}^3}$, with the manufactured solution
\begin{align}
  &\begin{aligned}
    &u_{\mathrm{ex}}\of{x} = \cos\of{k (x_1 - 3 x_2 + 2 x_3)}   \sin\of{k (1 + x_1)}   \\
    &\cdot \sin\of{k (1 - x_2)}   \sin\of{k (2 x_1 + x_2)}   \sin\of{k (3 x_1 - 2 x_2 + 2 x_3)} \eqdot
  \end{aligned}
      \intertext{The continuous right-hand side is set to}
      &f_{\mathrm{ex}} = \lambda u_{\mathrm{ex}} - \Delta u_{\mathrm{ex}}
\end{align}
and inhomogeneous~\dirichlet\ boundary conditions are imposed at all boundaries.
In the following, the parameter~$k$ is chosen as~${k = 5}$ and the~\helmholtz\ parameter as~${\lambda = 0}$, leading to the harder to solve \poisson's equation.

Four solvers are considered.
A conjugate gradient~(CG) solver using diagonal preconditioning, called~\solverdtcg~\cite{huismann_2017_condensation} which serves as an efficient baseline solver to compare to.
Then, a multigrid solver implementing~\prettyref{alg:mg}, called~\solvermg, using one pre- and one post-smoothing step, a~\krylov-accelerated version thereof based on~\prettyref{alg:cgmg} called~\solvermgcg, and, lastly, a \krylov-accelerated version with varying number of smoothing steps called~\solverkcvmg.

To test the scaling with the polynomial degree, the domain was discretized using~${\nelement=8^3}$ elements, with the polynomial degree being varied from~${\poly = 3}$ to~${\poly = 32}$.
As overresolution of the right-hand side might lead to faster convergence, the input data was instantiated using pseudo-random numbers, alleviating this effect.
The solvers were run 11 times, with the last 10 runs being averaged and the number of iterations to achieve a residual reduction by~$10^{-10}$, called $n_{10}$, measured.
From these, the convergence rate
\begin{align}
  \rho &= \sqrt[n_{10}]{\frac{\|\condr_{n_{10}}\|}{\|\condr_{0}\|}}
\end{align}
and the runtime per degree of freedom are computed.

\begin{figure}
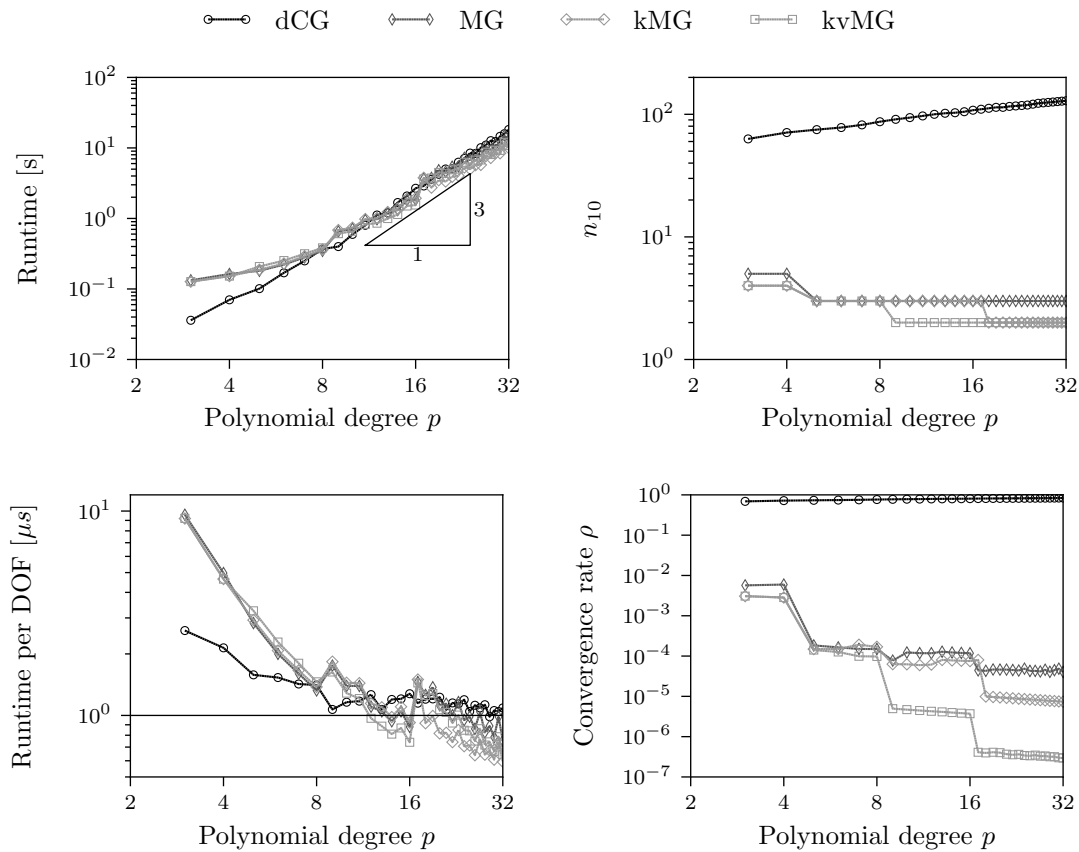

  \hspace*{\fill}
  \includegraphics{build/plot_runtime_test_helmholtz_PO_legend.pgf}
  \hspace*{\fill}\\
  \hspace*{\fill}
  \includegraphics{build/plot_runtime_test_helmholtz_PO_times.pgf}
  \hfill
  \includegraphics{build/plot_runtime_test_helmholtz_PO_iters.pgf}
  \hspace*{\fill}\\
  \hspace*{\fill}%
  \includegraphics{build/plot_runtime_test_helmholtz_PO_times_per_dof.pgf}
  \hfill
  \includegraphics{build/plot_runtime_test_helmholtz_PO_convergence.pgf}
  \hspace*{\fill}
  \caption{%
    Results for homogeneous meshes of~$\nelement = 8^3$ elements when varying the polynomial degree~$\poly$.
    Top left:~runtime of the solvers,
    top right:~number of iterations required to reduce the residual by 10 orders of magnitude,
    bottom left:~runtimes per degrees of freedom~(DOF),
    bottom right:~convergence rates of the solvers.}%
  \label{fig:homogeneous_poly_runtimes}
\end{figure}
\prettyref{fig:homogeneous_poly_runtimes} depicts the results.
While~the number of iterations increases for the CG solver, the multigrid variants exhibit a mostly constant iteration count.
The solver~\solvermg\ uses three iterations and has a slightly decreasing convergence rate of~$10^{-4}$, with the~\krylov\ acceleration improving matters and leading to only two iterations being used for high polynomial degrees.
This matches the performance of solvers with similar overlap, e.g.\,\cite{stiller_2017_multigrid}.
Applying the varying number of smoothing steps further improves the convergence rate, particularly when a new grid level is being introduced.
The consecutive addition of levels leads to the convergence rate decreasing from~$10^{-4}$ to~$10^{-5}$ and even to~$10^{-6}$ for~${\poly>16}$.
However, the better convergence rate does not directly translate to a lower runtime:
All multigrid solvers incur an overhead of a factor of four for low polynomial degrees.
At~$\poly = 8$, the runtime nearly equals that of the baseline solver~\solverdtcg.
However, due to the added level, the runtime increases again afterwards.
The multigrid solvers achieve a lower runtimes for~${\poly > 12}$, except for~${\poly > 17}$, where a new level is being introduced.
One has to bear in mind, that the number of elements is far lower than in practical computations, favoring the~CG solver.

\begin{figure}
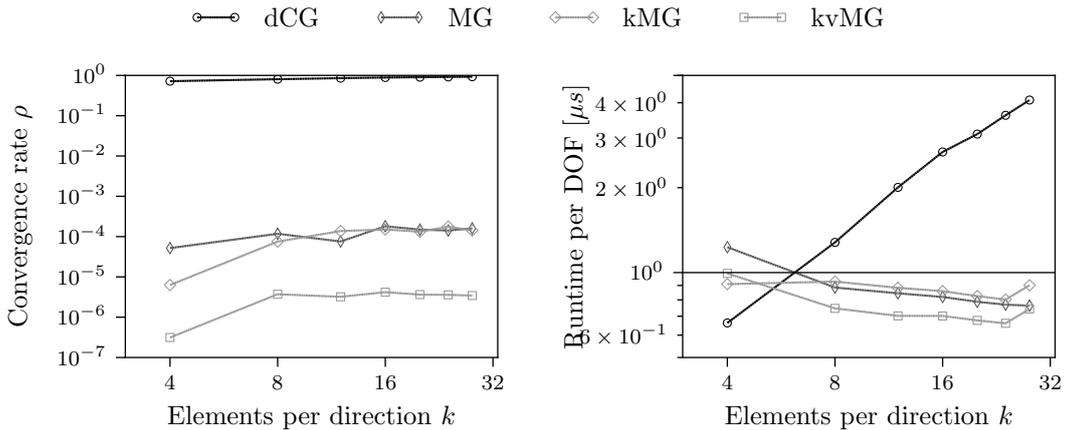

  \hspace*{\fill}
  \includegraphics{build/plot_runtime_test_helmholtz_NE_legend.pgf}
  \hspace*{\fill}\\
  \hspace*{\fill}
  \includegraphics{build/plot_runtime_test_helmholtz_NE_convergence.pgf}
  \hfill
  \includegraphics{build/plot_runtime_test_helmholtz_NE_times_per_dof.pgf}
  \hspace*{\fill}
  \caption{Convergence rates and runtimes per degree of freedom for the four solvers for homogeneous meshes of~$\nelement = k^3$ elements of polynomial degree~${\poly = 16}$.}
  \label{fig:homogeneous_nelement_runtimes}
\end{figure}

To validate that the multigrid solvers achieve a constant number of iterations when increasing the number of elements, the number of elements was varied at~${\poly=16}$.
\prettyref{fig:homogeneous_nelement_runtimes} depicts the convergence rate as well as the runtime per degree of freedom.
While the baseline solver~\solverdtcg\ yields an increasing number of iterations leading to an increasing runtime per degree of freedom, the multigrid solvers exhibit a nearly constant convergence rate, which translates to a constant iteration count.
The runtime, however, is not linear, it decreases until~$\nelement = 24^3$ and increases for~$\nelement=28^3$.
The first is an artifact of the vertex-based smoother: The ratio of vertices to elements decreases when increasing the elements per direction, leading to fewer evaluations per element and, hence, a lower runtime, whereas the latter stems from moving from using the RAM of one socket to non-uniform memory access across both sockets.

\subsection{Solver runtimes for anisotropic meshes}
So far, the solvers were only investigated for uniform meshes.
In simulations the resolution often needs to be adapted to the solution, leading to anisotropic or even stretched meshes.
To evaluate the effect of anisotropic meshes on the multigrid solvers the tests from~\cite{stiller_2017a_multigrid} were repeated: The aspect ratio~$AR$ of the mesh varies from~$AR = 1$ to~$AR = 48$, expanding the domain to
\begin{align}
  \domain = (0, 2 \pi \cdot AR) \times (0, \pi \left\lceil AR / 2 \right\rceil ) \times (0, 2 \pi) \eqdot
\end{align}
The domain was discretized using~${\nelement = 8^3}$ elements of polynomial degree~${\poly = 16}$ leading to a homogeneous mesh consisting of anisotropic brick-shaped elements.

\begin{figure}
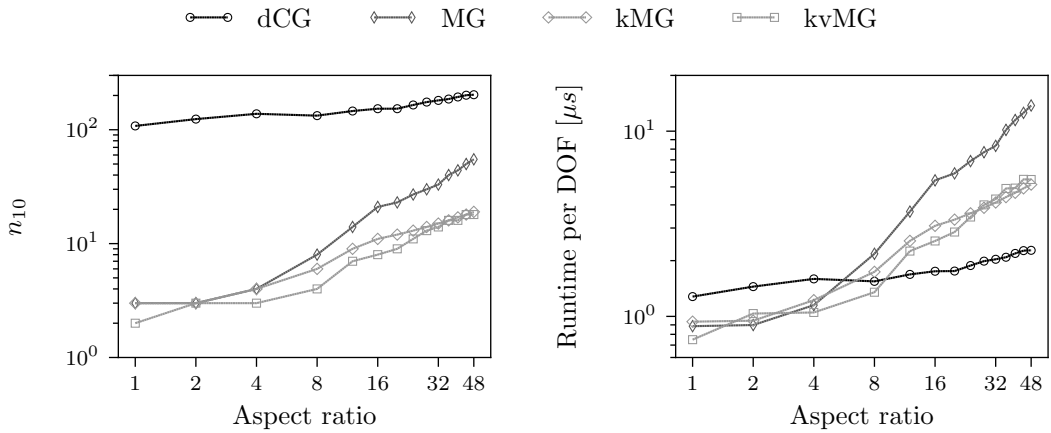

  \hspace*{\fill}
  \includegraphics{build/plot_runtime_test_helmholtz_AR_legend.pgf}
  \hspace*{\fill}\\
  \hspace*{\fill}
  \includegraphics{build/plot_runtime_test_helmholtz_AR_iters.pgf}
  \hfill
  \includegraphics{build/plot_runtime_test_helmholtz_AR_times_per_dof.pgf}
  \hspace*{\fill}
  \caption{Runtimes for the four solvers for anisitropic meshes of~$\nelement = 8^3$ elements of polynomial degree~${\poly = 16}$.}%
  \label{fig:anisotropic_runtimes}
\end{figure}

\prettyref{fig:anisotropic_runtimes} depicts the number of iterations and runtimes per degree of freedom of the solvers.
The locally-preconditioned solver~\solverdtcg\ exhibits a high robustness against increases in the aspect ratio, only requiring twice as long for an aspect ratio of~$AR=48$.
This is to be expected, as it bears similarity to so-called wirebasket solvers, even sharing their poly-logarithmic bounds regarding the polynomial degree~\cite{huismann_2017_wirebasket}.
The multigrid solvers do not fare as well.
For~\solvermg\ the number of iterations increases from three to sixty.
Applying \krylov-acceleration helps, but does not completely remove the effect.
While the solvers are very applicable for homogeneous meshes, their runtime increases rapidly for high-aspect ratios.
Increasing the number of pre- and post-smoothing cycles per level slightly mitigates the problem and keeps the number of iterations mostly stable until~$AR = 8$, but for higher aspect ratios a higher overlap is required~\cite{stiller_2017a_multigrid}.

\subsection{Solver runtimes for stretched meshes}
\begin{figure}
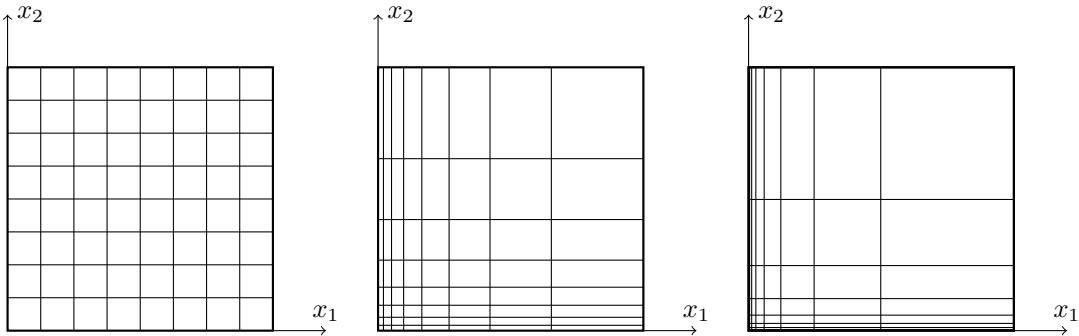

  \hspace*{\fill}
  \includegraphics[width=0.3 \textwidth]{images/non_homogeneous_grid_alpha_10.tikz}
  \hfill
  \includegraphics[width=0.3 \textwidth]{images/non_homogeneous_grid_alpha_15.tikz}
  \hfill
  \includegraphics[width=0.3 \textwidth]{images/non_homogeneous_grid_alpha_20.tikz}
  \hspace*{\fill}
  \caption{Cut through the~$x_3 = 0$ plane for the three meshes with constant expansion factor~$\alpha$.
    Left:~$\alpha = 1$, middle:~$\alpha = 1.5$, right:~${\alpha = 2}$}%
  \label{fig:stretched_meshes}
\end{figure}

%The last section considered homogeneous meshes.
%
In practice, applications often require a local mesh refinement, e.g.\,near walls, while still using isotropic elements in the flow.
To evaluate the robustness against varying the aspect ratio through the grid, the testcase from~\cite{huismann_2017_condensation} was adapted.
A grid consisting of~$8^3$ elements discretizes the domain~$\Omega = {(0,2\pi)}^3$ using a constant expansion factor~${\alpha \in \{1,1.5,2\}}$.
This leads to the three grids shown in~\prettyref{fig:stretched_meshes}.
Where~$\alpha = 1$, yields to a homogeneous mesh, $\alpha = 1.5$ stretches the grid in every cell, leading to a maximum aspect ratio of~${AR_{\mathrm{max}} = 17}$, while the last one is generated using~$\alpha = 2$ and has a maximum aspect ratio of~$AR_{\mathrm{max}} = 128$.
As frequently encountered in computational fluid dynamics, the grids are populated by a large variety of elements, ranging from boxfish, over plaice, to eels.
One has to keep in mind, that most solvers are not capable of handling such meshes at all for high polynomial degrees, as both, matrix-free operator evaluation, as well as matrix-free smoothers are required to handle the different element shapes.
Furthermore, the very large aspect ratios are detrimental to the performance of most solvers.

\prettyref{tab:stretched_poly_iterations} lists the number of iterations required to lower the residual by ten orders of magnitude.
The block-preconditioned solver~\solverdtcg\ exhibits only a slight increase in iterations when stretching the mesh from~$\alpha = 1$ to~$\alpha=1.5$ and~$\alpha=2$, using only fifty percent more iterations.
This ratio even gets slightly lowered when raising the polynomial degree.
For the multigrid solvers the situation is more complex.
For~$\alpha = 1$, the required number of iterations start at five or four, and get lowered to three or two, depending on the type of solver.
However, when raising the expansion factor to~$1.5$, the number of iteration increases fourfold for the traditional multigrid solver at low polynomial degrees and threefold for the \krylov-accelerated versions, with small improvements for large polynomial degrees.
For~${\alpha= 2}$ the effect gets more pronounced, with the multigrid solver~\solvermg\ requiring~$36$ iterations instead of~$5$ for low polynomial degrees, and~$12$ instead of~$3$ for~$\poly = 32$.
\krylov-acceleration lowers the number of iterations to~$15$ and~$8$, respectively.
However, increasing the number of smoothing steps is not beneficial anymore, again indicating that the overlap, not the amount of smoothing, is the main factor for attaining convergence.
\begin{table}
  \caption{Number of iterations required for reducing the residual by 10 orders of magnitude for the stretched grids using a constant expansion factor~$\alpha$.}%
  \label{tab:stretched_poly_iterations}
  \begin{center}
    \begin{tabular}{lrrrrrrrrrrrrrrr}\toprule
      && && \multicolumn{7}{c}{$\poly$}\\\cmidrule(lr){5-11}
      $\alpha$ &\phantom{a}& Solver &\phantom{a}& 4 &\phantom{a}& 8 &\phantom{a}& 16 &\phantom{a}& 32\\\midrule
      \input{build/table_alpha_niter}
      \bottomrule
    \end{tabular}
  \end{center}
\end{table}

%%% Local Variables:
%%% mode: latex
%%% TeX-master: "block_smoothers"
%%% End:

%% file: build/table_alpha_niter.tex
$1$&&dCG&& $71$ && $87$ && $108$ && $129$ \\
$1$&&MG&& $5$ && $3$ && $3$ && $3$ \\
$1$&&kMG&& $4$ && $3$ && $3$ && $2$ \\
$1$&&kvMG&& $4$ && $3$ && $2$ && $2$ \\\midrule
$1.5$&&dCG&& $98$ && $117$ && $126$ && $144$ \\
$1.5$&&MG&& $21$ && $11$ && $7$ && $5$ \\
$1.5$&&kMG&& $11$ && $8$ && $6$ && $4$ \\
$1.5$&&kvMG&& $11$ && $8$ && $5$ && $3$ \\\midrule
$2$&&dCG&& $105$ && $133$ && $158$ && $180$ \\
$2$&&MG&& $36$ && $26$ && $18$ && $12$ \\
$2$&&kMG&& $15$ && $13$ && $10$ && $8$ \\
$2$&&kvMG&& $15$ && $13$ && $10$ && $8$ \\

%% file: block_smoothers_extensions.tex
\section{Performance in simulations}\label{sec:discussion}

To solve the full~\navierstokes\ equations for large-scale problems many components are required, fast~\helmholtz\ solvers are just one, and all need to be parallelized.
To evaluate the performance in simulations, first the parallel performance of the multigrid solver is investigated then the performance in flow simulations are investigated computing the turbulent plane channel flow and the turbulent~\taylorgreen\ vortex.

\subsection{Parallelization}

To attain a good convergence rate, the presented solver combines a~$\poly$-multigrid approach with a vertex-based \schwarz-type smoother.
This kind of method requires data from the elements sharing the vertex.
A simple way to share that information in the parallel case are so-called ghost elements, which pad the subdomain of every process.
However, where the number of elements on a partition scales with~$\nelement = k^3$, the number of vertices scales with~$\nvertex = {(k+1)}^3$.
Furthermore, $k$ is small compared to traditional finite difference or finite volume methods.
For instance, $k=4$ results in a factor of four more vertices than elements and, hence, a factor of four in the work occurs -- in the worst case.
This leads to small subdomains being relatively expensive and thus renders a purely domain decomposition-based parallelization inefficient.
A two-level parallelization approach is required, with one coarse-grain level implementing the domain decomposition, and a fine-grain parallelization exploiting data parallelism inside the operators of a process, as done in~\cite{jin_2011_hybrid}.
In the present work, a domain decomposition layer is realized with MPI and combined with a second, fine-grain layer exploiting shared memory on CPUs via OpenMP\@.
Similarly, GPUs can be used in the second layer, e.g.\,using CUDA or OpenACC~\cite{kloeckner_warburton_nodal_dg_gpu}.

\begin{table}[]
  \centering
  \caption{Speedups and parallel efficiency for the four solver over the number of threads and the number of elements per direction on each process~$k$.}%
  \label{tab:speedup_test_one_node}
  \begin{tabular}{rrrrccccccccccc}\toprule
    &&&& \multicolumn{5}{c}{Speedup} && \multicolumn{5}{c}{Efficiency $[\%]$}\\\cmidrule(lr){5-9}\cmidrule(lr){11-15}
    &&&& \multicolumn{5}{c}{Number of threads} && \multicolumn{5}{c}{Number of threads}\\
    $k$ && Solver && 4&& 8 && 12 && 4 && 8 && 12 \\ \midrule
    \input{build/table_speedup.tex}\bottomrule
  \end{tabular}
\end{table}
To evaluate the efficiency of the parallelization, the test from section~\prettyref{sec:solver_runtimes} were repeated on one node using two~MPI processes, one per socket.
To allow for cuboidal subdomains with isotropic elements, the domain was set to~${\domain = (0,2 \cdot 2 \pi)\times {(0,2\pi)}^2}$ and decomposed in the~$x_1$-direction, with the number of elements scaled from~$2\cdot 8 \times 8^2$ to~$2 \cdot 16\times 16^2$ using~${\poly=16}$.
The number of threads per process and, hence, cores per process was varied from~$1$ to the maximum number of available cores per~CPU, which was~$12$ for the machine available, and the speedup over using only one thread per process determined.

\prettyref{tab:speedup_test_one_node} lists the resulting speedups and efficiencies.
For a small number of elements per process, the parallel efficiency declines rapidly.
The non-multigrid solver attains only~$60\ \%$ with eight cores and~$40\ \%$ with twelve, whereas the multigrid variants fare better with~$70\ \%$ and $55\ \%$, respectively.
Slight differences are present, with the~\krylov-accelerated versions being more efficient.
Increasing the number of elements for one process to~$12^3$ significantly increases the speedup and, hence, the efficiency which is~$80\ \%$ up to eight cores and~$65 \%$ for twelve.
Further increasing the number of elements per subdomain generates a higher speedup only for a small number of threads, indicating that the source of the low efficiency is not communication or latency.

The current implementation does generate acceptable but not optimal speedups.
However, one has to keep in mind that neither the L3 nor the memory bandwidth scales linearly with the active number of cores on the architecture utilized~\cite{hackenberg_2015_energy}.
The former scales up to a factor of ten and the latter to a factor of five compared to one core.
Furthermore, the current communication implementation represents a communication barrier.
Since not every~MPI library allows every thread to communicate, all but the master thread are idle during communication~\cite{mpi_standard}.
This can be remedied by overlapping computation and communication which allows to hide latency, albeit at a relatively high implementation cost~\cite{hindenlang_2012_dg}.

\subsection{Turbulent plane channel flow}

Until now the proposed multigrid solver has only been evaluated in analytical test cases.
While these allow to evaluate smoothing rates, they do not capture the behavior in flow simulations, where computing time, not the highest smoothing rate, is of paramount importance.
In flow solvers the computation of the pressure commonly takes the largest portion of the runtime, figures of it consuming up to $90\ \%$ are not unheard of~\cite{fehn_2018_sim}.
%
% As typically the runtime per degree of freedom increases with the polynomial degree, this inhibits the usage of high-order time-stepping schemes, which need multiple pressure solves per time step and become prohibitively expensive.
%
% Lowering the runtime of the pressure solver towards those of the explicit terms allows to use higher-order time-stepping techniques and, hence, attain the convergence order expected of the method.

The multigrid solver was implemented into the in-house flow solver~\specht.
The code is a testbed for heterogeneous parallelization as well as novel factorization techniques.
For time stepping it uses the incremental pressure-correction technique from~\cite{guermond_2006_projection}.
The method employs a backward-differencing scheme of second order accuracy and treats the convection terms explicitly via extrapolation, whereas the diffusion terms are treated implicitly.
One \poisson\ solve is required to project the intermediate velocity into the divergence-free space at the end of the time step.
Hence, one time step consists of evaluating the convection terms, solving a~\helmholtz\ equation for every velocity component, and then solving a \poisson\ equation for the velocity correction.

For spatial discretization structured grids of spectral elements using~\GLL\ points are employed.
The convection terms are evaluated using consistent integration, also referred to as overintegration, to eliminate aliasing problems.
In convection-dominated flows, the diffusion equations can be provided with a very good initial guess, allowing the usage of the baseline CG solver~\solverdtcg.
These, typically need less than ten iterations to solve, making them preferable to the multigrid methods.
The pressure solve, on the other hand, requires a multigrid approach for large-scale simulations and, here, the proposed multigrid method is utilized.

In this section the turbulent channel flow at~$Re_{\tau} = 180$ is considered~\cite{moser_1999_dns}.
The domain is set to~${\Omega = (0, 2\pi) \times (0, 2) \times (0, \pi)}$ and is periodic in~$x_1$ and~$x_3$ direction, whereas walls are present at~$x_2 = 0$ and~$x_2 = 2$.
A body force fixing the mean velocity to~$1$ drives the flow, with a PI controller computing the required force.
In combination with a bulk \reynolds\ number of~$5600$ and suitably disturbed initial conditions, this leads to a fully-developed turbulent flow near the~${Re_{\tau} = 180}$ mark, with~${Re_{\tau}}$ being a result of the simulations.

\begin{figure}[t]
  \centering
  \includegraphics{build/plot_channel_grid}\hfill \raisebox{.5cm}{\includegraphics[width=0.5\textwidth]{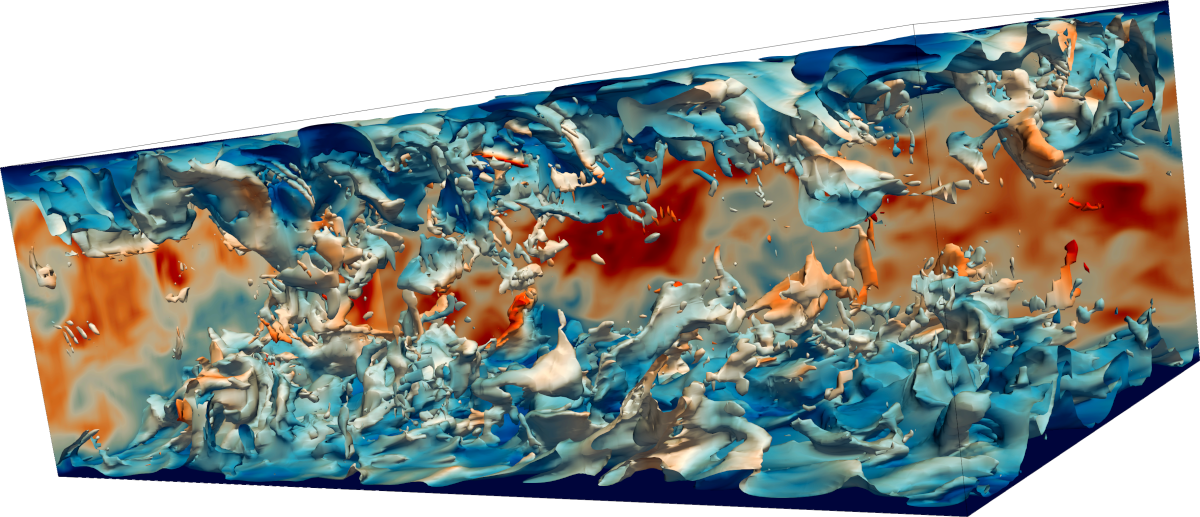}}
  \caption{Simulation of turbulent plane channel flow. Left: Cut through the~$x_1-x_2$-plane of the grid of spectral elements, only the element boundaries are shown. Right:~Isosurface of a transported passive scalar with value of~$1$ at the top and~$0$ at the bottom wall.}%
  \label{fig:grid_channel}
\end{figure}

The channel was discretized using~${16 \times 12 \times 6}$ spectral elements of degree~${\poly = 16}$.
To ensure that enough points are located in the boundary layer, the grid was generated with an expansion factor of~${\alpha= 1.2}$, refining it near the wall such that the first mesh point lies at~${y^{+}_{1} = 0.24}$ and the first eight points lie below~$y^{+} < 10$.
The resulting grid contains approximately~${4.7}$ million points, leading to~${19}$ million degrees of freedom, and is shown in combination with isosurfaces of a passive scalar in~\prettyref{fig:grid_channel}.
The simulation was performed on one node using two processes until a nearly steady state for the body force was reached.
To gain runtime data, the code was instrumented using~\scorep~\cite{knuepfer_2012_scorep} and run once using one thread per process and once using twelve threads.
After reaching a statistical steady state, $0.1$ dimensionless units in time were computed, requiring~$n_{\mathrm{time}}=429$ time steps with the~\helmholtz\ equations solved to a residual of~$10^{-10}$.

\begin{table}[t]
  \centering
  \caption{%
    Speedup for the plane channel flow test case.
    Setup data and total runtime obtaineed with the flow solver~\specht\ when using the new multigrid solver for computing~$0.1$ dimensionless units in time.
    The number of unknowns is computed as~$\ndof = 4 \poly^3 \nelement$.}%
  \label{tab:channel_times}
  \begin{tabular}{llrrr}\toprule
    Number of threads && 1 && 12\\\midrule
    Number of time steps~$n_{\mathrm{time}}$ && 429 && 429 \\
    Number of degrees of freedom~$\ndof$ && 18,874,368 && 18,874,368\\
    Number of cores~$n_{\mathrm{cores}}$ && 2 && 24 \\
    Runtime~$\twall$ [s] && 3336 && 491 \\
    $(\ndof \cdot n_{\mathrm{time}} ) / (\twall  \cdot n_{\mathrm{cores}})$ [1/s]  && 1,213,600 && 687,126\\\bottomrule
  \end{tabular}
\end{table}
\begin{table}[t]
  \centering
  \caption{%
    Accumulated runtimes for the time stepping procedure of~\specht\ when using the new multigrid solver for computing a time interval of~$0.1$ for the channel flow over the number of threads per process.
    Only components directly in the time-stepping procedure were profiled with~\scorep.}%
  \label{tab:channel_details}
  \begin{tabular}{llrrrrrrr}\toprule
                      &&  \multicolumn{3}{c}{1 thread} &&    \multicolumn{3}{c}{12 threads} \\ \cmidrule(lr){3-5} \cmidrule(lr){7-9}
    Component         && [s] && [\%] && [s] && [\%]\\\midrule
    Convection terms  && 2138 && 34.6  &&  2294 && 22.0 \\
    Diffusion solver  && 1454 && 23.6  &&  3160 && 30.3 \\
    \poisson\ solver  && 2407 && 39.0  &&  4551 && 43.6 \\
    Other             &&  175 &&  2.8  &&   395 &&  4.1\\\midrule
    Total             && 6174 && 100.0 && 10400 && 100.0\\\bottomrule
  \end{tabular}
\end{table}
\prettyref{tab:channel_times} summarizes the wall clock time and the runtime per degree of freedom, while~\prettyref{tab:channel_details} lists the contributions of different components of the flow solver for two processes when decomposing the domain along the~$x_1$-direction.
For one thread per process the runtime consists mainly of three contributions:
The diffusion step, which only takes a quarter of the runtime, the convection terms, which take a third, and the computation of the pressure, requiring nearly~$40\ \%$ of the runtime.
For this case the mission is accomplished: Treating the pressure is as cheap as treating the convection terms.
The ratio would even be better for homogeneous grids, as here the pressure solver still needs only three iterations.
In these simulations, the code reaches a throughput of over~$1,200,000$ degrees of freedom per core and second.
When using twelve cores, the convection terms parallelize very well, as they require very few memory accesses, leading to nearly the same time accumulated in these routines.
However, this does not hold for the implicit solvers where the parallel efficiency is less than~$50\ \%$ for the diffusion solves and a bit above it for the multigrid solver.
This leads to an increase of the percentage in CPU time required by the pressure solver to~$44\ \%$, which is twice the value for the convection terms.
As a result, the throughput per core drops to~$680,000$ degrees of freedom per core and second.

\subsection{Turbulent~\taylorgreen\ vortex benchmark}
\begin{figure}[h]
  \includegraphics[width=0.3\textwidth]{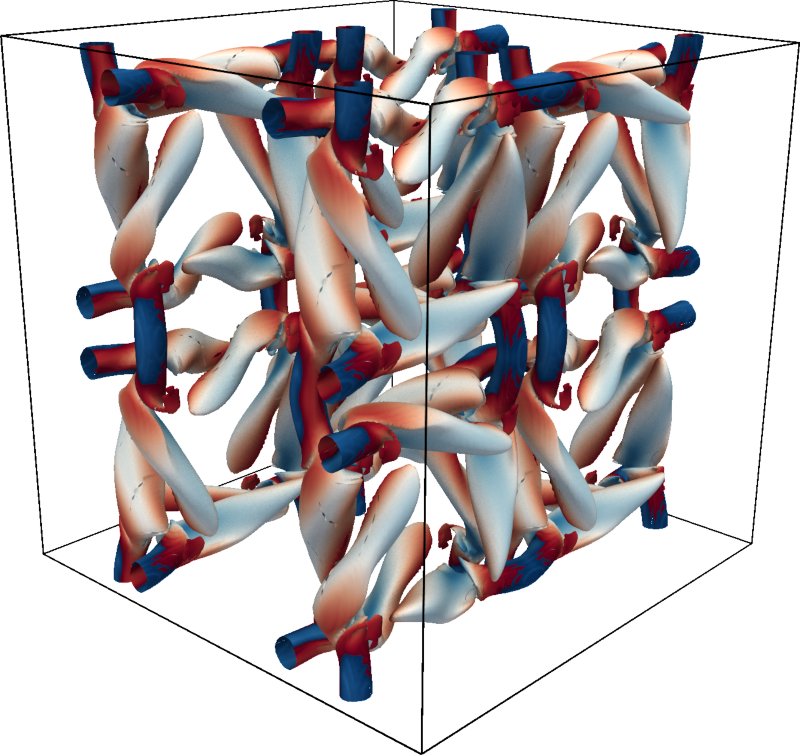}
  \hfill
  \includegraphics[width=0.3\textwidth]{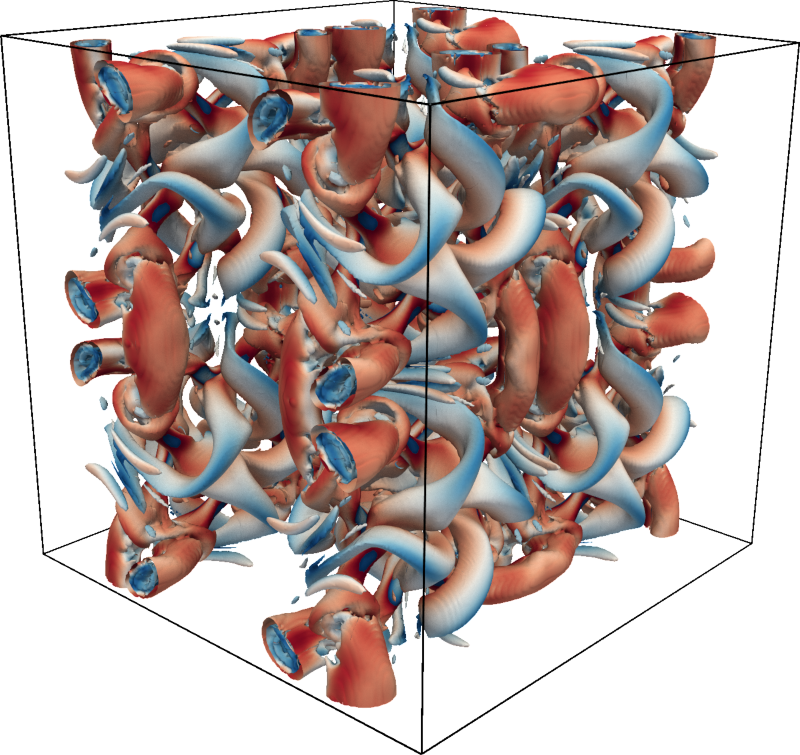}
  \hfill
  \includegraphics[width=0.3\textwidth]{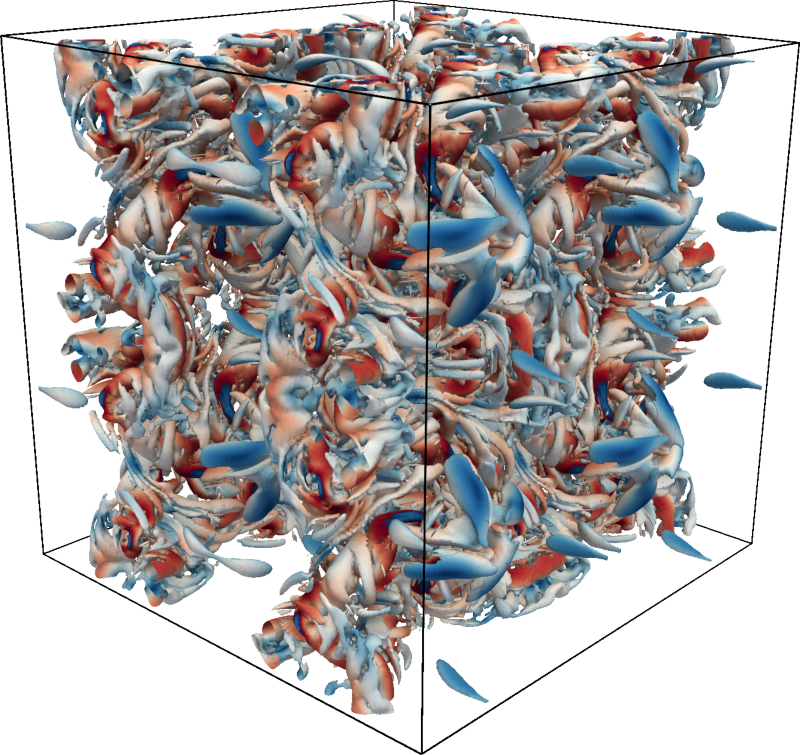}
  \caption{%
    Isosurfaces of the~$\lambda_2$ vortex criterion for the~\taylorgreen\ vortex at~${Re = 1600}$ with~${\lambda_2 = -1.5}$.
    The isosurfaces are colored with the magnitude of the velocity vector, where blue corresponds zero and red to a magnitude of~$U_0$.
    The data was taken from a simulation using~${\nelement=16^3}$ spectral elements of polynomial degree~${\poly=16}$.
    Left:~${t = 5  T_0}$, middle:~${t = 7  T_0}$, right:~${t = 9  T_0}$.}%
  \label{fig:turbulent_taylor_green}
\end{figure}

To further evaluate the efficiency of the code the underresolved turbulent \taylorgreen\ vortex benchmark is considered~\cite{gassner_2013_sim, fehn_2018_sim}.
For a length scale~$L$, a reference velocity~$U_0$, and a periodic domain~${\Omega = (-L\pi,L\pi)}^3$ the initial conditions for the velocity components are
\begin{subequations}
  \begin{align}
    u_1 \of{\vec{x}} &= +U_0 \sin\of{\frac{x_1}{L}} \cos\of{\frac{x_2}{L}} \cos\of{\frac{x_3}{L}}\\
    u_2 \of{\vec{x}} &= -U_0 \cos\of{\frac{x_1}{L}} \sin\of{\frac{x_2}{L}} \cos\of{\frac{x_3}{L}}\\
    u_3 \of{\vec{x}} &= 0 \eqdot
  \end{align}
\end{subequations}
At a~\reynolds\ number of~$Re = U_0 L / \nu = 1600$, where~$\nu$ is the kinematic viscosity, the flow is highly unstable and quickly transitions to turbulence~\cite{green_1937_fluidmech}, as shown in~\prettyref{fig:turbulent_taylor_green}.

Four homogeneous meshes are considered.
First, a grid using~$\nelement = 16^3$ elements of polynomial degree~${\poly=8}$.
The second one contains the same number of degrees of freedom but with~$\nelement = 8^3$ and~${\poly=16}$.
The third grid is finer grid with~$\nelement=16^3$ and~$\poly=8$.
The fourth grid has~$\nelement=16^3$ at~${\poly=16}$.
The simulations were carried out until a simulation time of~${T=20 L / U_0 = 20 T_0}$.
A constant time step was imposed and set according to the~spectral element CFL condition~\cite{karniadakis_1999_sem}
\begin{align}
  \Delta t &= C_{\mathrm{CFL}} \frac{\max_{i,e}|h_{i,e}|}{\poly^2} \max _{\Omega} | \vec{u}| \eqdot
\end{align}
Here, a CFL number of~$C_{\mathrm{CFL}} = 0.125$ was imposed and the~\helmholtz\ equations are solved to an absolute tolerance of~$10^{-10}$ for the residual.

In this benchmark, the grids were deliberately chosen to be very coarse and, so that they are not capable of capturing all features of the flow.
In the discontinuous~\galerkin\ methods utilized in~\cite{gassner_2013_sim, fehn_2018_sim}, the flux formulation leads to an implicit subgrid-scale~(SGS) model which generates the required dissipation and stabilization.
The present work, however, uses continuous elements with less numerical dissipation and no inherent~SGS model.
As a remedy, the spectral vanishing viscosity model~(SVV) was employed.
It modifies the~\laplace\ matrix of the velocities, adding more dissipation for high polynomial degrees, resulting in a certain stabilization~\cite{koal_2012_svv}.
The power kernel by~\moura~\cite{moura_2016_svv} was chosen, with the model parameters set to~$\poly_{\mathrm{SVV}} = \poly/2$ and~$\varepsilon_{\mathrm{SVV}} = 0.01$.
It confines the viscous effects of the SVV to higher polynomial modes and leads to improved accuracy for high polynomial degrees.

\begin{table}[t]
  \caption{%
    Grids utilized for the turbulent~\taylorgreen\ benchmark in conjunction with the respective number of degrees of freedom~$\ndof$, number of data points~$\ndofstar$, number of time steps~$\ntime$, wall clock time~$\twall$, number of cores~$\ncores$, CPU time, and computational throughput~${(\ntime \cdot \ndofstar) / (\twall \cdot \ncores)}$.}%
  \label{tab:taylor_green_times}
  \begin{centering}
  \begin{tabular}{lrrrrrrrrrrrrrrrr}\toprule
    % $\poly$                && $8$ && $8$ && $16$ && $16$ && $16$   \\
    && \multicolumn{3}{c}{$\poly=8$} && \multicolumn{5}{c}{$\poly=16$}   \\\cmidrule(lr){3-5} \cmidrule(lr){7-11}
    $\nelement$            && $16^3$ && $32^3$ && $4^3$ && $8^3$ && $16^3$  \\\midrule
    $\ndof$                && 8,388,608  && 67,108,864 && 1,048,576 &&  8,388,608 && 67,108,864  \\
    $\ndofstar$            && 11,943,936 && 95,551,488 && 1,257,728 && 11,943,936 && 80,494,592  \\
    $\ntime$               && 26,076     && 52,152 && 26,076 && 52,152 && 104,304  \\\midrule
    $e^2_{E_k} \cdot 10^3$ &&  3.57 && 0.487 && 17.8 && 1.37 && 0.229  \\
    Runtime $\twall$ [\SI{}{\second}]   && 13,964 && 74,104 &&  2,297 && 25,855 && 110,058  \\
    $\ncores$              && 24 && 96 && 24 && 24 && 96   \\
    CPU time [CPUh]        && 93 && 1976 && 15 && 172 && 2935   \\
    $(\ndofstar \cdot n_{\mathrm{time}} ) / (\twall  \cdot n_{\mathrm{cores}})$ [1/s]            && 929,356 && 700,480 && 594,827 && 845,664 && 794,650  \\\bottomrule
  \end{tabular}
  \end{centering}
\end{table}

\begin{figure}[t]
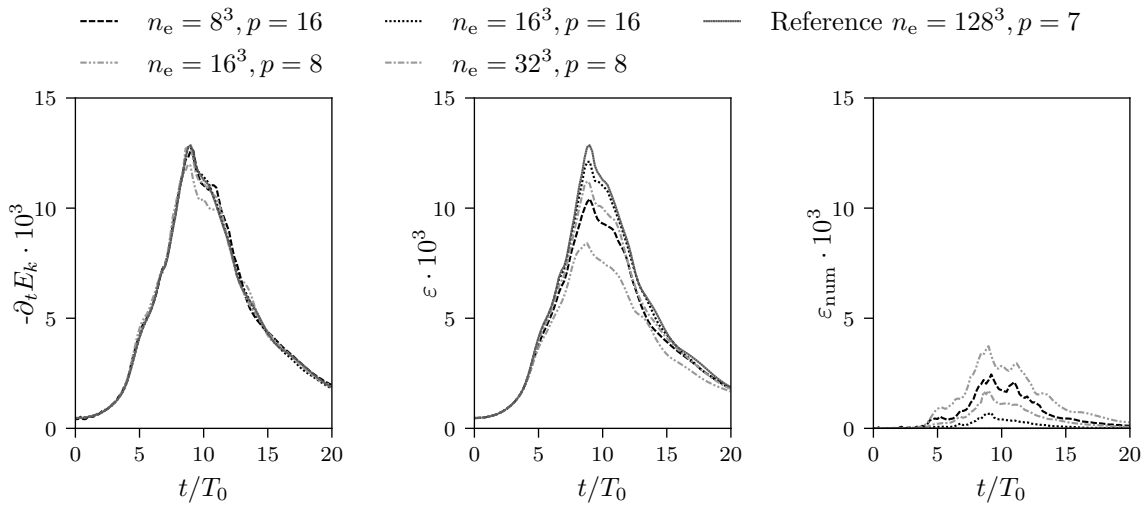

  \hspace*{\fill}
  \includegraphics{build/plot_taylor_green_legend.pgf}
  \hspace*{\fill}\\
  \includegraphics{build/plot_taylor_green_energy.pgf}
  \hfill
  \includegraphics{build/plot_taylor_green_dissipation.pgf}
  \hfill
  \includegraphics{build/plot_taylor_green_numerical_dissipation.pgf}
  \caption{%
    Results for the turbulent~\taylorgreen\ vortex. Left: Time derivative of the mean kinetic energy, middle:~mean dissipation rate captured by the grid, right:~numerical dissipation.
    Reference data courtesy of~M.\,Kronbichler~\cite{fehn_2018_sim}.}%
  \label{fig:taylor_green_energy}
\end{figure}

\prettyref{fig:taylor_green_energy} depicts the derivative of the mean kinetic energy in the subdomain~$E_k$ over time combined with the respective dissipation rate~$\varepsilon$ and their difference, the numerical dissipation.
The results are compared to DNS data from~\cite{fehn_2018_sim}.
The coarsest grid with${\nelement = 16^3}$ and~${\poly = 8}$ is capable of initially resolving the flow, but starting from~${t = 4 T_0}$ deviations are present and the peak in energy loss is not obtained correctly.
The deviations stem from smaller molecular dissipation and larger numerical dissipation, which peaks at a third of the reference dissipation rate.
Increasing the number of elements and, hence, the number of degrees of freedom, leads to more of the dissipation rate being resolved.
However, keeping the number of grid points constant and increasing the polynomial degree instead decreases the error more noticeably.

To quantify the accuracy of the results, the relative~$L_2$-error of the time derivative of the kinetic energy~$E_k$ is computed as
\begin{align}
  e^2_{E_k} &= \frac{\int\limits_0^T{\left(\partial_{t} E_k\of{\tau} - \partial_{t} E_{k,\mathrm{ref}}\of{\tau}\right)}^2\mathrm{d}\tau}{\int\limits_0^T{\left(\partial_{t} E_{k,\mathrm{ref}}\of{\tau}\right)}^2\mathrm{d}\tau} \eqcomma
\end{align}
where~${E_{k,\mathrm{ref}}}$ is the kinetic energy from the reference data.
\prettyref{tab:taylor_green_times} lists the accuracy of the simulations in combination with the number of time steps, the number of degrees of freedom, the error in~$E_k$, and the achieved computational throughput.
To attain comparability with~\cite{fehn_2018_sim}, the throughput is computed from the number of element-local grid points, ${(\poly+1)}^3 \nelement$, times four for the number of solution variables.
The discontinuous formulation converges towards the result from the continuous formulation, as the solution is continuous, negating the extra degrees of freedom allowed for in the discontinuous case.
For a constant number of time steps and, hence, the same time step width, using more degrees of freedom leads to a smaller error.
This indicates that the testcase is spatially underresolved as the error does not only depend on the time step width.
Except for the coarsest mesh using~$4^3$ elements on~$24$ cores the computational throughput is~$800,000$ timestep times the number of data points per second and core.
This is observed for both, $\poly =8$ as well as~$\poly=16$.
In all simulations the multigrid solver used for the pressure equation required only two iterations.
As a result, the higher computational cost for the convection terms at higher polynomial degrees are offset by a more efficient multigrid cycle.
Furthermore, slight parallelization losses are present when increasing the number of nodes from one to four, leading to a higher throughput for computations with fewer number of elements.
Compared to~\cite{fehn_2018_sim}, where the testcase was run on comparable hardware, a factor of two in throughput was achieved here.
To the knowledge of the authors, this makes~\specht\ the fastest solver for incompressible flow employing high polynomial degrees, at the time of writing.

% \subsection{Extension to discontinuous~Galerkin methods}
% \begin{itemize}
% \item In~\cite{stiller_2016_multigrid}, a tensor-product version of the~\helmholtz\ operator for~DG was presented.
% \item With the same arguments, the fast diagonalization can be combined with static condensation in the~DG formulation.
% \item leads to a linearly scaling operator
% \item albeit with double the number of DOF per subdomain
% \item[$\Rightarrow$] operator costs twice as much
% \item requires static condensation in the elements
% \item possible using the techniques from~\cite{huismann_2017_condensation}
% \item however, the operator is also twice as expensive as the fluxes require the derivative of the variable on the element boundary

% \end{itemize}
%%% Local Variables:
%%% mode: latex
%%% TeX-master: "block_smoothers"
%%% End:

%% file: build/table_speedup.tex
8&&dCG&& $3.18$ && $4.51$ && $4.60$ && $ 79$ && $ 56$ && $ 38$ \\
8&&MG&& $3.25$ && $4.85$ && $5.85$ && $ 81$ && $ 61$ && $ 49$ \\
8&&kMG&& $3.47$ && $5.40$ && $6.26$ && $ 87$ && $ 68$ && $ 52$ \\
8&&kvMG&& $3.42$ && $5.69$ && $6.58$ && $ 85$ && $ 71$ && $ 55$ \\\midrule
12&&dCG&& $3.60$ && $4.95$ && $5.31$ && $ 90$ && $ 62$ && $ 44$ \\
12&&MG&& $3.62$ && $6.31$ && $7.84$ && $ 91$ && $ 79$ && $ 65$ \\
12&&kMG&& $3.61$ && $6.27$ && $7.76$ && $ 90$ && $ 78$ && $ 65$ \\
12&&kvMG&& $3.62$ && $6.24$ && $7.83$ && $ 91$ && $ 78$ && $ 65$ \\\midrule
16&&dCG&& $3.66$ && $4.95$ && $5.20$ && $ 91$ && $ 62$ && $ 43$ \\
16&&MG&& $3.72$ && $6.45$ && $7.94$ && $ 93$ && $ 81$ && $ 66$ \\
16&&kMG&& $3.70$ && $6.38$ && $7.78$ && $ 93$ && $ 80$ && $ 65$ \\
16&&kvMG&& $3.70$ && $6.38$ && $7.84$ && $ 92$ && $ 80$ && $ 65$ \\

%% file: block_smoothers_conclusion.tex
\section{Conclusions}\label{sec:conclusions}
This paper presented a novel factorization for the inverse of the~\helmholtz\ operator on the~$2 \times 2 \times 2$ element block.
The statically condensed system was embedded into the full one, allowing to use the fast diagonalization as matrix-free inverse and factorize it, reducing the operator complexity from~$\order{\poly^4}$ to~$\order{\poly^3}$.
Hence, a linearly scaling inverse for the~$2^3$ element block was derived.
Then, runtime tests were conducted proving that the inverse indeed scales linearly and performs better than matrix-matrix multiplication with the inverse, even outpacing the matrix-matrix version starting from~${\poly=5}$ and beating the fast diagonalization for the full system for every relevant polynomial degree~$\poly$.
This was generalized to larger blocks in~\prettyref{sec:element_centered}.

Using the linearly scaling inverse as main building block for an overlapping \schwarz\ smoother, a $p$-multigrid solver utilizing static condensation was proposed.
The multigrid cycle scales linearly when combining the residual evaluation from~\cite{huismann_2017_condensation} with the linearly scaling inverse, only leaving pre- and postprocessing for static condensation with super-linear contributions.

Tests for the solver were conducted. To reduce the residual by a factor~$10^{-10}$, it required less than four iterations, in most cases two to three.
The linearly scaling operators lead to very high efficiency for the multigrid cycle, so that for the solver uses less than one microsecond per unknown over a wide parameter range when computing on one core.
Furthermore, the runtime spent in pre- and postprocessing is relatively small, so that the solver scales linearly with the number of degrees of freedom in the range of polynomial degrees tested and is expected to do so until~${\poly = 48}$.

Parallelization studies were conducted and the performance in full~\navierstokes\ simulations was evaluated.
Where traditional incompressible flow solvers require up to~$90\ \%$ of the runtime to compute the pressure, the new multigrid solver lowers that margin to the portion required to evaluate the explicitly treated convection terms, albeit with slight losses generated by the parallelization.
This makes the code~\specht\ faster than any other present high order incompressible flow solver.

Future work could be dedicated to expanding the multigrid solver towards the discontinuous \galerkin\ method, increasing its range of applicability.
Furthermore, efficiency gains are to be expected from fusing operators and using cache-blocking to decrease the impact of the limited memory bandwidth.
Lastly, a GPU implementation for the finest grid could substantially boost the performance.

%%% Local Variables:
%%% mode: latex
%%% TeX-master: "block_smoothers"
%%% End:

%% file: block_smoothers.bbl
\begin{thebibliography}{10}

\bibitem{atak_2016_sim}
M.~Atak, A.~Beck, T.~Bolemann, D.~Flad, H.~Frank, and C.-D. Munz.
\newblock High fidelity scale-resolving computational fluid dynamics using the
  high order discontinuous {Galerkin} spectral element method.
\newblock In {\em High Performance Computing in Science and Engineering{
  \'{}}15}, pages 511--530. Springer, 2016.

\bibitem{beck_2014_dg}
A.~D. Beck, T.~Bolemann, D.~Flad, H.~Frank, G.~J. Gassner, F.~Hindenlang, and
  C.-D. Munz.
\newblock High-order discontinuous {Galerkin} spectral element methods for
  transitional and turbulent flow simulations.
\newblock {\em International Journal for Numerical Methods in Fluids},
  76(8):522--548, 2014.

\bibitem{bramble_1995_mg}
J.~Bramble.
\newblock {\em Multigrid methods}.
\newblock Pitman Res. Notes Math. Ser. 294. Longman Scientific \& Technical,
  Harlow, UK, 1995.

\bibitem{brandt_1982_mg}
A.~Brandt.
\newblock Guide to multigrid development.
\newblock In {\em Multigrid Methods}, volume 960 of {\em Lecture Notes in
  Mathematics}, pages 220--312. Springer Berlin/Heidelberg, 1982.

\bibitem{couzy_1995_condensation}
W.~Couzy and M.~Deville.
\newblock A fast {Schur} complement method for the spectral element
  discretization of the incompressible {Navier}-{Stokes} equations.
\newblock {\em Journal of Computational Physics}, 116(1):135 -- 142, 1995.

\bibitem{deville_2002_sem}
M.~Deville, P.~Fischer, and E.~Mund.
\newblock {\em High-Order Methods for Incompressible Fluid Flow}.
\newblock Cambridge University Press, 2002.

\bibitem{dryja_1994_substructuring}
M.~Dryja, B.~F. Smith, and O.~B. Widlund.
\newblock {Schwarz} analysis of iterative substructuring algorithms for
  elliptic problems in three dimensions.
\newblock {\em SIAM Journal on Numerical Analysis}, 31(6):1662--1694, 1994.

\bibitem{fehn_2018_sim}
N.~Fehn, W.~A. Wall, and M.~Kronbichler.
\newblock Efficiency of high-performance discontinuous {Galerkin} spectral
  element methods for under-resolved turbulent incompressible flows.
\newblock {\em International Journal for Numerical Methods in Fluids}, 2018.

\bibitem{feng_2001_mg}
X.~Feng and O.~A. Karakashian.
\newblock Two-level additive {Schwarz} methods for a discontinuous {Galerkin}
  approximation of second order elliptic problems.
\newblock {\em SIAM Journal on Numerical Analysis}, 39(4):1343--1365, 2001.

\bibitem{gander_2008_schwarz}
M.~J. Gander et~al.
\newblock {Schwarz} methods over the course of time.
\newblock {\em Electronic Transactions on Numerical Analysis}, 31(5):228--255,
  2008.

\bibitem{gassner_2013_sim}
G.~J. Gassner and A.~D. Beck.
\newblock On the accuracy of high-order discretizations for underresolved
  turbulence simulations.
\newblock {\em Theoretical and Computational Fluid Dynamics}, 27(3-4):221--237,
  2013.

\bibitem{golub_1999_ipcg}
G.~H. Golub and Q.~Ye.
\newblock Inexact preconditioned conjugate gradient method with inner-outer
  iteration.
\newblock {\em SIAM Journal on Scientific Computing}, 21(4):1305--1320, 1999.

\bibitem{green_1937_fluidmech}
A.~E. Green and G.~I. Taylor.
\newblock Mechanism of the production of small eddies from larger ones.
\newblock {\em Proceedings of the Royal Society of London A}, 158, 1937.

\bibitem{guermond_2006_projection}
J.-L. Guermond, P.~Minev, and J.~Shen.
\newblock An overview of projection methods for incompressible flows.
\newblock {\em Computer methods in applied mechanics and engineering},
  195(44):6011--6045, 2006.

\bibitem{hackbusch_1985_mg}
W.~Hackbusch.
\newblock {\em Multigrid Methods and Applications}, volume~4 of {\em
  Computational Mathematics}.
\newblock Springer, 1985.

\bibitem{hackenberg_2015_energy}
D.~Hackenberg, R.~Sch{\"o}ne, T.~Ilsche, D.~Molka, J.~Schuchart, and R.~Geyer.
\newblock An energy efficiency feature survey of the {Intel} {Haswell}
  processor.
\newblock In {\em Parallel Distributed Processing Symposium Workshops (IPDPSW),
  2015 IEEE International}, 2015.

\bibitem{haupt_2017_mg}
L.~Haupt.
\newblock {\em Erweiterte mathematische Methoden zur Simulation von turbulenten
  Strömungsvorgängen auf parallelen Rechnern}.
\newblock PhD thesis, Centre for Information Services and High Performance
  Computing (ZIH), TU Dresden, Dresden, 2017.
\newblock (in German).

\bibitem{haupt_2013_mg}
L.~Haupt, J.~Stiller, and W.~E. Nagel.
\newblock A fast spectral element solver combining static condensation and
  multigrid techniques.
\newblock {\em Journal of Computational Physics}, 255(0):384 -- 395, 2013.

\bibitem{hestenes_1952_cg}
M.~R. Hestenes and E.~Stiefel.
\newblock Methods of conjugate gradients for solving linear systems.
\newblock {\em Journal of Research of the National Bureau of Standards},
  49(6):409--436, 1952.

\bibitem{hindenlang_2012_dg}
F.~Hindenlang, G.~J. Gassner, C.~Altmann, A.~Beck, M.~Staudenmaier, and C.-D.
  Munz.
\newblock Explicit discontinuous {Galerkin} methods for unsteady problems.
\newblock {\em Computers \& Fluids}, 61:86--93, 2012.

\bibitem{huismann_2017_wirebasket}
I.~Huismann, J.~Stiller, and J.~Fr{\"o}hlich.
\newblock Building blocks for a leading edge high-order flow solver.
\newblock {\em Proceedings in Applied Mathematics and Mechanics},
  17(1):129--132, 2017.

\bibitem{huismann_2017_condensation}
I.~Huismann, J.~Stiller, and J.~Fr{\"o}hlich.
\newblock Factorizing the factorization {\textendash} a spectral-element solver
  for elliptic equations with linear operation count.
\newblock {\em Journal of Computational Physics}, 346:437--448, 2017.

\bibitem{jin_2011_hybrid}
H.~Jin, D.~Jespersen, P.~Mehrotra, R.~Biswas, L.~Huang, and B.~Chapman.
\newblock High performance computing using {MPI} and {OpenMP} on multi-core
  parallel systems.
\newblock {\em Parallel Computing}, 37:562--575, 2011.

\bibitem{karniadakis_1999_sem}
G.~Karniadakis and S.~Sherwin.
\newblock {\em Spectral/hp Element Methods for {CFD}}.
\newblock Oxford University Press, 1999.

\bibitem{kloeckner_warburton_nodal_dg_gpu}
A.~Klöckner, T.~Warburton, J.~Bridge, and J.~Hesthaven.
\newblock Nodal discontinuous {Galerkin} methods on graphics processors.
\newblock {\em Journal of Computational Physics}, 228(21):7863 -- 7882, 2009.

\bibitem{knuepfer_2012_scorep}
A.~Kn{\"u}pfer, C.~R{\"o}ssel, D.~a. Mey, S.~Biersdorff, K.~Diethelm,
  D.~Eschweiler, M.~Geimer, M.~Gerndt, D.~Lorenz, A.~Malony, W.~E. Nagel,
  Y.~Oleynik, P.~Philippen, P.~Saviankou, D.~Schmidl, S.~Shende,
  R.~Tsch{\"u}ter, M.~Wagner, B.~Wesarg, and F.~Wolf.
\newblock Score-{P}: A joint performance measurement run-time infrastructure
  for {Periscope}, {Scalasca}, {TAU}, and {Vampir}.
\newblock In H.~Brunst, M.~S. M{\"u}ller, W.~E. Nagel, and M.~M. Resch,
  editors, {\em Tools for High Performance Computing 2011}, pages 79--91,
  Berlin, Heidelberg, 2012. Springer Berlin Heidelberg.

\bibitem{koal_2012_svv}
K.~Koal, J.~Stiller, and H.~M. Blackburn.
\newblock Adapting the spectral vanishing viscosity method for large-eddy
  simulations in cylindrical configurations.
\newblock {\em Journal of Computational Physics}, 231(8):3389--3405, 2012.

\bibitem{lombard_2015_sim}
J.-E.~W. Lombard, D.~Moxey, S.~J. Sherwin, J.~F.~A. Hoessler, S.~Dhandapani,
  and M.~J. Taylor.
\newblock Implicit large-eddy simulation of a wingtip vortex.
\newblock {\em AIAA Journal}, pages 1 -- 13, 2015.

\bibitem{lottes_2005_mg}
J.~W. Lottes and P.~F. Fischer.
\newblock Hybrid multigrid/{Schwarz} algorithms for the spectral element
  method.
\newblock {\em Journal of Scientific Computing}, 24(1):45--78, 2005.

\bibitem{lynch_1964_tensors}
R.~Lynch, J.~Rice, and D.~Thomas.
\newblock Direct solution of partial difference equations by tensor product
  methods.
\newblock {\em Numerische Mathematik}, 6(1):185--199, 1964.

\bibitem{merzari_2013_sim}
E.~Merzari, W.~Pointer, and P.~Fischer.
\newblock Numerical simulation and proper orthogonal decomposition of the flow
  in a counter-flow t-junction.
\newblock {\em Journal of Fluids Engineering}, 135(9):091304, 2013.

\bibitem{moser_1999_dns}
R.~D. Moser, J.~Kim, and N.~N. Mansour.
\newblock Direct numerical simulation of turbulent channel flow up to
  ${Re_{\tau}= 590}$.
\newblock {\em Physics of Fluids}, 11(4):943--945, 1999.

\bibitem{moura_2016_svv}
R.~Moura, S.~Sherwin, and J.~Peir{\'o}.
\newblock Eigensolution analysis of spectral/hp continuous galerkin
  approximations to advection--diffusion problems: Insights into spectral
  vanishing viscosity.
\newblock {\em Journal of Computational Physics}, 307:401--422, 2016.

\bibitem{oosterlee_1998_mg}
C.~W. Oosterlee and T.~Washio.
\newblock An evaluation of parallel multigrid as a solver and a preconditioner
  for singularly perturbed problems.
\newblock {\em SIAM Journal on Scientific Computing}, 19(1):87--110, 1998.

\bibitem{pasquetti_2009_multigrid}
R.~Pasquetti and F.~Rapetti.
\newblock p-multigrid method for {Fekete}-{Gauss} spectral element
  approximations of elliptic problems.
\newblock {\em Communications in Computational Physics}, 5(2-4):667--682, Feb
  2009.

\bibitem{patera_1984_sem}
A.~T. Patera.
\newblock A spectral element method for fluid dynamics: laminar flow in a
  channel expansion.
\newblock {\em Journal of Computational Physics}, 54(3):468 -- 488, 1984.

\bibitem{ronquist_1987_mg}
E.~M. R{\o}nquist and A.~T. Patera.
\newblock Spectral element multigrid. {I}. formulation and numerical results.
\newblock {\em Journal of Scientific Computing}, 2(4):389--406, 1987.

\bibitem{serson_2017_cfd}
D.~Serson, J.~R. Meneghini, and S.~J. Sherwin.
\newblock Direct numerical simulations of the flow around wings with spanwise
  waviness.
\newblock {\em Journal of Fluid Mechanics}, 826:714--731, 2017.

\bibitem{sherwin_2001_sem}
S.~J. Sherwin and M.~Casarin.
\newblock Low-energy basis preconditioning for elliptic substructured solvers
  based on unstructured spectral/hp element discretization.
\newblock {\em Journal of Computational Physics}, 171(1):394--417, 2001.

\bibitem{stiller_2016_multigrid}
J.~Stiller.
\newblock Robust multigrid for high-order discontinuous {Galerkin} methods: {A}
  fast {Poisson} solver suitable for high-aspect ratio {Cartesian} grids.
\newblock {\em Journal of Computational Physics}, 327:317--336, 2016.

\bibitem{stiller_2017_multigrid}
J.~Stiller.
\newblock Nonuniformly weighted {Schwarz} smoothers for spectral element
  multigrid.
\newblock {\em Journal of Scientific Computing}, 72(1):81--96, 2017.

\bibitem{stiller_2017a_multigrid}
J.~Stiller.
\newblock Robust multigrid for cartesian interior penalty {DG} formulations of
  the {Poisson} equation in 3d.
\newblock In {\em Spectral and High Order Methods for Partial Differential
  Equations ICOSAHOM 2016}, pages 189--201. Springer, 2017.

\bibitem{mpi_standard}
{The {MPI} Forum}.
\newblock {MPI}: A message passing interface version 3.0, 2012.

\bibitem{woodward_1998_newton}
C.~S. Woodward.
\newblock A {Newton}-{Krylov}-multigrid solver for variably saturated flow
  problems.
\newblock {\em WIT Transactions on Ecology and the Environment}, 24, 1998.

\end{thebibliography}
